\documentclass[a4paper,11pt]{article}

\usepackage[english]{babel}
\usepackage[latin1]{inputenc}

\usepackage{lmodern}
\usepackage{textcomp}
\usepackage{enumerate}
\usepackage{graphicx,epsfig,psfrag}
\usepackage{mathrsfs}
\usepackage{amsmath,amsfonts,amssymb}
\usepackage{color}
\usepackage{xcolor}
\usepackage{fancybox}
\usepackage{pgf,tikz}
\usetikzlibrary{arrows} 
\usepackage{animate}
\definecolor{rltred}{rgb}{0.75,0,0} 
	\definecolor{rltgreen}{rgb}{0,0.5,0}
	\definecolor{oneblue}{rgb}{0,0,0.75}
	\definecolor{marron}{rgb}{0.64,0.16,0.16}
	\definecolor{forestgreen}{rgb}{0.13,0.54,0.13}
	\definecolor{purple}{rgb}{0.62,0.12,0.94}
	\definecolor{dockerblue}{rgb}{0.11,0.56,0.98}
	\definecolor{freeblue}{rgb}{0.25,0.41,0.88}
	\definecolor{myblue}{rgb}{0,0.2,0.4}
	\definecolor{ccqqtt}{rgb}{0.8,0,0.2}

\newtheorem{lem}{Lemma}[section]
\newtheorem{coro}{Corollary}[section]
\newtheorem{prop}{Proposition}[section]
\newtheorem{defin}{Definition}[section]
\newtheorem{theo}{Theorem}[section]

\newenvironment{dem}{\vskip 2mm\noindent \textbf{Proof:}}{\hfill $\square$ \vskip 2mm \noindent}

\title{Existence and percolation results for stopped germ-grain models with unbounded velocities}

\newcommand{\Leb}{\textrm{Leb}}
\newcommand{\1}{\ensuremath{\mbox{\rm 1\kern-0.23em I}}}

\renewcommand{\H}{\text{Hex}}

\begin{document}

\author{
David Coupier$^{1}$, David Dereudre$^{2}$
 and Simon Le Stum$^{3}$ \\
{\normalsize{\em }}
 }

\maketitle

\footnotetext[1]{\;Universit\'e Polytechnique des Hauts-de-France,  David.Coupier@uphf.fr }
 \footnotetext[2]{\;University of Lille, david.dereudre@univ-lille.fr }
 \footnotetext[3]{\;University of Lille, simon.lestum@univ-lille.fr }

\begin{abstract}
We investigate the existence and first percolation properties of  general stopped germ-grain models. They are defined via a random set of germs generated by a homogeneous planar Poisson point process in  $\mathbf{R}^{2}$. From each germ, a grain, composed by a random number of branches, grows. This grain stops growing whenever one of its branches hits another grain. The classical and historical example is the  line segment model for which the grains are segments growing in a random direction in $ [0,2\pi)$ with random velocity. In the bilateral line segment model the segments grow in both directions. Other examples are considered here such as the Brownian model where the  branches are simply given by independent  Brownian motions in  $\mathbf{R}^{2}$.
The existence of such dynamics for an infinite number of germs is not obvious and our first result ensures it in a very general setting. In particular the existence of the line segment model is proved as soon as the random velocity admits a moment of order 4 which extends the result by Daley et al (Theorem 4.3 in \cite{daley2014two}) for bounded velocity. Our result covers also the Brownian dynamic model.
In the second part of the paper, we show that the line segment model with random velocity admitting a super exponential moment does not percolate. This improves a recent result (Theorem 3.2 \cite{coupier2016absence}) in the case of bounded velocity.

    \bigskip


\end{abstract}

%

\section{Introduction}
\label{SectionIntro}

Consider a stationary Poisson point process $\mathbf{X}$ in $\mathbf{R}^{2}$ where each germ $\xi\in\mathbf{X}$ is marked by a uniformly random direction $\Theta$ on $[0,2\pi]$. At time $0$, each germ gives off a growing line segment at unit rate, in the associated direction. One of the ends of the growing line segment is $\xi$, and the other one determines the stop of the segment when it hits another line segment. The almost sure existence of a unique stopped system of non-overlapping finite line segments has been proved by Daley et al (Theorem 4.3 in \cite{daley2014two}). The authors conjectured the absence of percolation for this model and this has been recently solved by Coupier et al. (Theorem 3.2 in \cite{coupier2016absence}). Percolation means here the existence of an unbounded connected component produced by the union of stopped segments. When the random velocity $\mathbf{V}$ takes its values in a compact set, the existence and absence of percolation come from  slight modifications   of theorems mentioned above. In the case of random unbounded velocities $\mathbf{V}$, the problem is much more complicated because very quick segments may destroy the locality of the dynamic. As a consequence of our main theorems, we establish the existence  of the stopped line segment model as soon as $\mathbf{E}(\mathbf{V}^{4})<+\infty$ and the absence of percolation if there exists $s>1$ such that $\mathbf{E}(e^{\mathbf{V}^{s}})<+\infty$. Obviously this super-exponential moment condition is quite restrictive but it covers the Gaussian random velocities case  corresponding to the classical Maxwellian velocity distribution for an ideal gas in statistical mechanics.

In a general setting, a planar germ-grain dynamic is a system of growing particles defined on a Poisson point process. The set of germs is distributed by a standard Poisson point process in $\mathbf{R}^{2}$. A typical grain is defined as a finite collection of random processes identically distributed in $\mathbf{R}^{2}$. Each grain can be viewed as a finite collection of growing branches coming from the associated germ. Any growing grain ceases its propagation when one of its extremities hits another grain. Given an initial infinite configuration of growing grains, the existence and uniqueness of such a stopped germ-grain dynamic is not guaranteed. In Theorem \ref{theo:ex} we establish the almost sure existence and uniqueness under a very general assumption involving the fourth moment of the expansion of the grain (see Assumption \eqref{equ:suffi}). As mentioned above this assumption is equivalent to a finite moment of order four (i.e. $\mathbf{E}(\mathbf{V}^{4})<+\infty$) in the line segment model. Our result covers also the existence of the Brownian model where each grain is a collection of independent Brownian motions in $\mathbf{R}^{2}$. 

Let us now turn to our results on the absence of percolation for some stopped germ-grain models. As mentioned above we prove such absence of percolation for the line segment model in the case $\mathbf{E}(e^{\mathbf{V}^s})<+\infty$ for some $s>1$.  For simplicity we deal only with the line segment model although the result could be easily extended to several stopped germ-grain models having similar geometric properties. "Similar geometric properties" means that the direction of expansion of the grain is random but does change during the evolution. The randomness is completely  encoded in the initial condition. For instance our results do not cover the case where the grains grow with unpredictable directions as in the Brownian model. The absence of percolation for the Brownian model is still a conjecture today.
 
The absence of  percolation for the line segment model with unbounded velocities is an extension of a recent result in \cite{coupier2016absence} with bounded velocities. This non-trivial extension required the development of new concepts and ideas which we explain briefly now. From any stopped germ-grain model, we associate an outdegree-one graph where the vertices are simply the points of the Poisson point process itself and the outgoing edges are defined by pointing out the stopping germ. With this formalism, the absence of percolation for a stopped germ-grain model is equivalent to the absence of percolation for its associated outdegree-one graph, for which a general theorem is developed in \cite{coupier2016absence} under the so-called assumptions Loop and Shield.  Roughly speaking, the Loop assumption means that any forward branch merges on a loop provided that the Poisson point process is augmented with a finite collection of well-chosen points. This assumption is still true in our model with unbounded velocity. In our setting, the shield assumption means that a large square box $[-n,n]^2$ has high probability (i.e. the probability tends to $1$ when $n\to \infty$) of not being crossed by a segment from the left to the right, uniformly with respect to the configuration outside the box. Obviously this property does not occur for unbounded velocities since it is always possible to build a very quick segment starting outside the box and crossing any line structure inside the box before its formation. So the major issue here is to show that germs with high velocities do not pollute so much the space $\mathbb R^2$ in order to apply the general strategy developed  in \cite{coupier2016absence} in the non-polluted part of the space. In particular, we need that this non-polluted part is large enough and percolates. The super exponential moment condition is then required to build shield blocks with high probability in the non-polluted domain.

The paper is organized as follows. In Section 2, we provide a precise description of stopped germ-grain models and we give the three main examples. In Section 3 we present our results on existence and absence of percolation. Sections 4 and 5 are devoted to the proofs.

\section{The stopped germ-grain model}
\label{sec:germ}

This section is devoted to the notion of stopped germ-grain model. The set of \textbf{germs} is generated by a homogeneous planar Poisson point process. From each germ, a \textbf{grain} is growing, made up with a random number of \textbf{branches} which are identically distributed but not necessarily independent. The grains are assumed to be independent from each other. Any grain ceases to grow whenever one of its branches hits another grain. 

Our first example of stopped germ-grain model is the famous line segment model introduced in \cite{daley2014two} in which each grain simply corresponds to one (unilateral or bilateral) segment growing with a constant velocity.

\subsection{The germ-grain model}

All models in this paper take place in the Euclidean space $\mathbf{R}^{2}$. The associated Lebesgue measure is denoted by $\Leb$. The intensity of the (homogeneous) Poisson point process generating the germs is $\lambda \Leb$ with $\lambda>0$. The numbers of branches per grain are i.i.d. positive random integers with distribution $\delta$ and whose generic r.v. is denoted by $K$. The state space for the collection of grains is
$$
\mathbf{F} := \left(\mathscr{C}(\mathbf{R}^{+},\mathbf{R}^{2})\right)^{\mathbf{N}} = \left\{ (f_{n})_{n\ge 0}\ ;\ \forall n\ge 0,\ f_{n}:\mathbf{R}^{+} \to \mathbf{R}^{2}\text{ is continuous}\right\} ~.
$$
{The space $\mathbf{F}$  allocates  an infinite number of continuous paths for each grain although only a finite number (but random) of paths is useful. This random number is furnished by another mark (see below). This formalism is a bit over technical but that it is the way we find to set rigorously the space of marks.}

The product space $\mathbf{F}$ is equipped with the cylindric $\sigma$-algebra $\mathcal{S}$ where each set $\mathscr{C}(\mathbf{R}^{+},\mathbf{R}^{2})$ is equipped with the Borel $\sigma$-algebra generated by the uniform convergence on compact sets. Let us note that we consider a infinite collection of branches for each grain although we will only used a finite number of them. 

Let us consider a probability measure $\mathcal{L}$ on the measurable space $(\mathbf{F},\mathcal{S})$ satisfying two assumptions:
\begin{itemize}
\item Given $Y=(Y_{0},Y_{1},\dots,Y_{n},\dots)$ distributed according to $\mathcal{L}$, all marginals are identically distributed (but not necessarily independent); 
\item For any index $i$, a.s. $Y_{i}(0)=0$, i.e. all the branches of a given grain start their trajectories from the corresponding germ.
\end{itemize}

The general mark space $\mathcal{M}$ of our model is defined by $\mathcal{M}=\mathbf{N}^{*}\times\mathbf{F}$, where $\mathbf{N}^{*}$ denotes the set of positive integers, and the configuration space $\mathcal{C}^{\mathcal{M}}$ on $\mathbf{R}^{2}$ with marks in $\mathcal{M}$ is defined by
$$
\mathcal{C}^{\mathcal{M}} = \left\{\varphi\subset\mathbf{R}^{2}\times\mathcal{M}\ ;\ \#(\varphi\cap(\Lambda\times\mathcal{M}))<\infty,\ \text{for any bounded }\Lambda\subset\mathbf{R}^{2}\right\} ~.
$$
It is equipped with the $\sigma$-algebra $\mathscr{F}$ generated by the counting events $P_{(A,n)}=\lbrace \varphi\in\mathcal{C}^{\mathcal{M}}\ ;\ \#(\varphi\cap A)\le n\rbrace$ for all $n\ge 0$ and $A$ in the sigma field $\mathcal{P}(\mathbf{N}^{*})\otimes \mathcal{S}$.

{A $(\lambda,\delta,\mathcal{L})$-germ-grain model can be roughly described as follows. From any germ (or Poisson point) $\xi$, a random number $K_{\xi}$ (with distribution $\delta$) of grains are started, each of them with distribution $\mathcal{L}$. In a more formal way:}

\begin{defin}
\label{def:GG}
A $(\lambda,\delta,\mathcal{L})$-\textbf{germ-grain model} is a Poisson point process on $\mathcal{C}^{\mathcal{M}}$ with intensity $\lambda\Leb\otimes\delta\otimes\mathcal{L}$.
\end{defin}

Let $\mathbf{X}$ be such a germ-grain model. The associated germ process is denoted by
$$
\mathbf{X}_{\text{germs}} := \left\{\xi : \, (\xi,\cdot,\cdot) \in \mathbf{X}\right\} \subset \mathbf{R}^2 ~.
$$
Given a marked point $(\xi,k,Y)\in\mathbf{X}$, the integer $k$ corresponds to the number of branches starting at the germ $\xi$. And, in the countable collection of random variables $(Y_{i})_{i\ge 0}$, only the $k$ first paths are used in the dynamic. For a time $t\ge 0$ and a marked point $x=(\xi,k,Y)$, we define the corresponding \textbf{grain} until time $t$ by
$$
\text{Grain}(x,t) := \left\{ \xi+Y_{i}(s) : \, 0\le i\le k-1 , \ s\in[0,t) \right\} ~.
$$
Also, the \textbf{extremity} at time $t$ of the previous grain is
$$
\text{H}(x,t) := \left\{ \xi+Y_{i}(t) : \, 0\le i\le k-1\right\}
$$
whose elements are called \textbf{particles}. Let us note that a grain does not contain its extremity. This convention has been chosen in order to simplify the stopping rules in Definition \ref{def:stopcong}.

\subsection{Examples of germ-grain model}
\label{sec:examples}

\subsubsection{Unilateral line segment model (Model 1)}
\label{sect:ULS}

This model is directly inspired by recent works \cite{coupier2016absence,daley2014two} on planar line segment dynamics. First, each grain is made up with only one branch, i.e. the distribution $\delta$ is the Dirac measure on $1$. To specify the probability measure $\mathcal{L}$, we need two independent r.v.'s. On the one hand, $\Theta$ which is uniformly distributed on $[0,2\pi]$ gives the direction of the line segment. On the other hand, $V$ on $\mathbf{R}^{*}_{+}$ is the growth velocity of the line segment. Then the variable $Y_1$ is defined by:
$$
\forall t \ge0 , \, Y_{1}(t) := t V \big( \cos\Theta , \sin\Theta \big)
$$
and the  probability measure $\mathcal{L}$ is obtained as the law of the sequence $Y:=(Y_{1},Y_{1},\dots)$. In other words, the germ-grain $(\xi,1,Y)$ is a single and unilateral line segment growing from the germ $\xi$ according to the direction $\Theta$ and with velocity $V$. 

For simplicity, in several places of the paper, the unilateral line segment model is described by a Poisson point process $\mathbf{X}$ on $\mathbf{R}^{2}\times[0,2\pi]\times\mathbf{R}^{*}_{+}$ with intensity $\lambda\Leb\otimes\Xi\otimes\mu$, where $\Xi$ is the uniform distribution on $[0,2\pi]$ and $\mu$ is the law of $V$. We will then use the notation
$$
\mathbf{X} = \bigcup_{\xi\in\mathbf{X}_{\text{germs}}} \left(\xi,\Theta_{\xi},V_{\xi} \right) ~.
$$

\subsubsection{Bilateral line segment model (Model 2)}

The bilateral line segment model was introduced in \cite{daley2014two} as well. Its definition also involves the r.v.'s $\Theta$, $V$ and $Y_{1}$ from the previous section. The probability measure $\mathcal{L}$ is now the law of the sequence
$Y:=(Y_{1},-Y_{1},Y_{1},Y_{1},\dots)$ and the distribution $\delta$ is the Dirac measure on $2$. Hence, the germ-grain $(\xi,2,Y)$ is a bilateral line segment growing from $\xi$ according to the two opposite directions $\Theta$ and $-\Theta$ with the same velocity $V$.

This example highlights the possibility of having dependent marginals in $\mathcal{L}$. Several other examples could be constructed similarly.

\subsubsection{Brownian model (Model 3)}

In this model, any given grain is made up with a random number $K$ of branches which are driven by i.i.d. Brownian trajectories. For technical reasons, we assume that $\mathbf{E}(K^{2})$ is finite. Thus let us set $B:=(B_{i})_{i\ge 0}$ a sequence of i.i.d. Brownian motions in $\mathbf{R}^{2}$ starting from $0$. The probability measure $\mathcal{L}$ is then simply the distribution of $B$. Hence, a germ-grain is completely determined by the triplet $(\xi,K,B)$.  

To illustrate the variety of possible models considered here, we could also imagine a single Brownian path starting from each germ, or exactly two paths but one reflected (w.r.t. the germ) from the other, etc.


\subsection{Stopped germ-grain model}

This section is devoted to the notion of {\bf stopped germ-grain model}. Instead of a dynamical definition as described in the Introduction, we prefer here to define a stopped germ-grain model through the concept of lifetime. This point of view was already used to define the Lilypond model in \cite{daley2005descending}.

\begin{defin}
\label{def:stopcong}
Let $\varphi\in\mathcal{C}^{\mathcal{M}}$ be a configuration. An \textbf{exploration} is a function defined on $\varphi$ associating a lifetime to each marked point.
$$
f : \varphi \longrightarrow (0,+\infty] ~.
$$
An exploration $f$ of $\varphi$ is said \textbf{stopped} if it satisfies the two following conditions:
\begin{itemize}
\item[$(i)$] \textbf{(Hardcore property).} $\forall x\not= y\in\varphi$ then $\text{Grain}(x,f(x))\cap\text{Grain}(y,f(y)) = \emptyset$;
\item[$(ii)$] \textbf{(Uniqueness of the stopping grain).} $\forall x\in\varphi$ with $f(x)<+\infty$ then  $\exists ! y\in\varphi\setminus\!\{x\}$ such that $H(x,f(x))\cap\text{Grain}(y,f(y)) \not= \emptyset$. 
\end{itemize}
\end{defin}

Item $(ii)$ asserts that either a marked point $x\in\varphi$ has an infinite lifetime, i.e. $f(x)=+\infty$, and will be never stopped or it will be eventually stopped, i.e. $f(x)<+\infty$, but by only one grain.

\begin{defin}
\label{def:pog}
A $(\lambda,\delta,\mathscr{L})$-germ-grain model $\mathbf{X}$ is said \textbf{stopped} if it a.s. admits a unique \textbf{stopped exploration}, denoted by $f_{\mathbf{X}}$. Moreover, $\mathbf{X}$ satisfies the \textbf{finite time property} if
\begin{equation}
\label{equa:poiss}
\mathbf{P} \big( \forall x\in\mathbf{X} , \, f_{\mathbf{X}}(x) < +\infty \big) = 1 ~.
\end{equation}
\end{defin}

The existence and uniqueness of a stopped exploration is not obvious in general because of the infinite number of marked points. We will prove in Corollary \ref{coro:models} and Proposition \ref{prop:POG} that the three germ-grain models introduced in Section \ref{sec:examples} are stopped and satisfy the finite time property.

It is worth pointing out here that the finite time property is crucial in our work since it allows us to interpret any stopped germ-grain model with such a property as a \textbf{Poisson outdegree-one graph} (see Section \ref{sect:Step1}). Indeed, each marked point has a finite lifetime and then admits a unique outgoing vertex corresponding to its stopping grain.

\section{Results}

\subsection{Existence of stopped germ-grain models}
\label{sec:ex}

In this section we state a general result (Theorem \ref{theo:ex}) ensuring the existence of stopped germ-grain models. Our main assumption involves the modulus of continuity for branches of a given grain. Precisely, for $\mathbf{X}$ a $(\lambda,\delta,\mathcal{L})$-germ-grain model where $K$ is a $\delta$-distributed r.v. and $Y=(Y_{i})_{i\ge 0}$ is a $\mathcal{L}$-distributed sequence of paths, and for all $t,t'\ge 0$, we set:
\begin{equation}
\label{MaxRadius}
M_{t,t'} := \max_{0\le k\le\mathbf{K}-1} \, \sup_{0\le s\le t'} \, \| Y_{k}(t+s) - Y_{k}(t) \| ~.
\end{equation}

{Roughly speeaking, the quantity $M_{t,t'}$ is the random variable of the maximal displacement, during the time interval $[t,t+t']$, of all branches of a typical grain with mark $(K,Y)$.}

\begin{theo}
\label{theo:ex}
Let $\mathbf{X}$ be a $(\lambda,\delta,\mathcal{L})$-germ-grain model such that $\mathbf{E}_{\delta}(K^{2})<+\infty$ and
\begin{equation}
\label{equ:suffi}
\lim_{t'\rightarrow 0}\ \sup_{t\ge 0}\mathbf{E}\left(M_{t,t'}^{4}\right)=0.
\end{equation}
Suppose also that, for any $U_{1},U_{2}$ two independent random paths on $\mathscr{C}(\mathbf{R}^{+},\mathbf{R}^{2})$ distributed as the marginals of $\mathcal{L}$,
\begin{equation}
\label{assregular}
\mbox{a.s. } \, \Leb\left(\{U_{1}(t),\ t\ge 0\}\right) = 0 \; \mbox{ and } \; \Leb\left(\{U_{1}(t)-U_{2}(t),\ t\ge 0\}\right) = 0 ~.
\end{equation} 
Then $\mathbf{X}$ is a stopped germ-grain model.
\end{theo}

The Assumption \eqref{assregular} is technical and satisfied for most natural models. It actually ensures that pathological random paths do not occur. Theorem \ref{theo:ex} will be proved in Section \ref{Section:existence}.

As a consequence of Theorem \ref{theo:ex}, we easily get that the three germ-grain models presented in Section \ref{sec:examples} are stopped.

\begin{coro}
\label{coro:models}
The unilateral and bilateral line segment models (Models  1 and 2) with $\mathbf{E}(V^{4})<+\infty$ and the Brownian model (Model 3) are stopped germ-grain models.
\end{coro}

Let us note that  D.J. Daley et al have established in \cite{daley2014two} the existence of the bilateral line segment model with constant velocity. They proved the finite time property as well.

\begin{dem}
For both line segment models, we have $M_{t,t'}=V t'$ and Assumption (\ref{equ:suffi}) is obviously satisfied. The Assumption \eqref{assregular} is easily checked as well.

Let us focus on the Brownian model. By stationarity of increments for Brownian paths, Assumption (\ref{equ:suffi}) is equivalent to
\begin{equation}
\label{eq:mike}
\lim_{t'\rightarrow 0} \mathbf{E} \left( M_{0,t'}^{4} \right) = 0 ~.
\end{equation}
Recall that
$$
M_{0,t'} = \max_{0\le k\le K-1} \sup_{0\le s\le t'} \|B_{k}(s)\| ~,
$$
where $(B_{k})_{0\le k\le K-1}$ is a collection of $K$ independent Brownian motions. Using the scaling property, we can write $\mathbf{E}\big( M_{0,t'}^{4} \big) = (t')^{2} \mathbf{E}\big( M_{0,1}^{4}\big)$. It then remains to prove that $\mathbf{E}\big( M_{0,1}^{4}\big)$ is finite. This simply comes from
$$
\mathbf{E}\left( M_{0,1}^{4} \right) = \mathbf{E} \left(\max_{0\le k\le K-1} W_{k}^{4} \right)\le \mathbf{E} \left(\sum_{k=0}^{K-1} W_{k}^{4} \right) = \mathbf{E}(K) \mathbf{E}(W_0^4) < +\infty ~,
$$
{where $W_{k}:=\sup_{0\le s\le 1} \|B_{k}(s)\|$ admits a finite fourth moment.} Finally, the Brownian model also satisfied Assumption \eqref{assregular} since the Lebesgue measure of the bi-dimensional Brownian path in $\mathbb{R}^2$ is a.s. equal to $0$.
\end{dem}

\subsection{Absence of percolation}

Let us focus on the graph produced by a stopped germ-grain model. 

\begin{defin}
\label{def:percolation}
For any $(\lambda,\delta,\mathcal{L})$ germ-grain model $\mathbf{X}$, we set
$$
\Sigma(\mathbf{X}) := \bigcup_{x\in\mathbf{X}} \text{Grain}(x,f_{\mathbf{X}}(x)) ~.
$$
We say that $\mathbf{X}$ percolates if $\Sigma(\mathbf{X})$ contains an unbounded connected component. 
\end{defin}

We conjecture the absence of percolation for a large class of stopped germ-grain models. This conjecture is supported by the underlying outdegree-one graph structure of stopped germ-grain models with finite time property (see Section \ref{sect:Step1}). The absence of percolation has been conjectured in \cite{daley2014two} and proved in \cite{coupier2016absence} for the line segment model with constant velocity. Our initial ambition was to state the absence of percolation for a large class of stopped germ-grain models but we did not succeed in this task and the conjecture is still largely open today. Nevertheless, we significantly improve the result of \cite{coupier2016absence} by allowing velocities of line segments to be random and especially unbounded.

Before giving our main theorem below (Theorem \ref{theo:perco}), let us briefly discuss the finite time property given in Definition \ref{def:pog}. The percolation question is relevant only if the studied stopped germ-grain model satisfies the finite time property. Otherwise, some of its grains have an infinite lifetime and then are (generally) unbounded. Besides, a large class of stopped germ-grain models should satisfy the finite time property (but we have not investigated such a general result in the present paper). This is the case for the three models introduced in Section \ref{sec:examples}.

\begin{prop}
\label{prop:POG}
The unilateral and bilateral line segment models (Models  1 and 2) with $\mathbf{E}(V^{4})<+\infty$ and the Brownian model (Model 3) are stopped germ-grain models satisfying the finite time property.
\end{prop}

The proof of Proposition \ref{prop:POG} is given at the end of this section. Here is our second main result stating the absence of percolation for the unilateral line-segment model with unbounded velocities. The proof is given in Section \ref{Section_perco}.

\begin{theo}
\label{theo:perco}
Assume that there exists $s>1$ such that $\mathbf{E}(\exp(V^{s}))<+\infty$. Then the unilateral line segment model does not percolate.
\end{theo}

\begin{dem}{\textbf{(of Proposition \ref{prop:POG})}.}
First, let us focus on the unilateral line segment model which is described by a Poisson point process $\mathbf{X}$ on $\mathbf{R}^{2}\times[0,2\pi]\times\mathbf{R}^{*}_{+}$ where $\mathbf{R}^{*}_{+}$ is the set of positive real numbers. Let $\gamma$ be a typical point located at the origin: $\gamma=(0,\Theta,V)$ where $\Theta$ is a uniform r.v. on $[0,2\pi]$ and $\mathbf{E}(V^{4})<+\infty$. Let us proceed by contradiction, assuming that with positive probability the typical point $\gamma$ is never stopped: $\mathbf{P}\big(f_{\mathbf{X}\cup\lbrace\gamma\rbrace}(\gamma)=+\infty\big)>0$. The Campbell-Mecke formula implies that, for any Borel set $\Delta\subset\mathbf{R}^{2}$,
\begin{equation}
\label{eq:Mecke}
\mathbf{E}\left(\sum_{x\in\mathbf{X}_{\Delta}} \1_{\lbrace f_{\mathbf{X}}(x)=+\infty \rbrace}\right)=\lambda\Leb(\Delta)\mathbf{P}\left(f_{\mathbf{X}\cup\lbrace\gamma\rbrace}(\gamma)=+\infty\right) \, > \, 0 ~,
\end{equation}
where $\mathbf{X}_{\Delta}=\{ x=(\xi,\Theta,V)\in\mathbf{X} : \xi\in\Delta\}$. For an introduction to the Palm theory, the reader may refer to \cite{chiu2013stochastic}. The isotropy of the r.v. $\Theta$ and (\ref{eq:Mecke}) allow us to state that for any $n\ge 1$
$$
\mathbf{P}\left(\exists x=(\xi,\Theta,V)\in\mathbf{X}\ ;\ f_{\mathbf{X}}(x)=\infty,\; \xi\in[n,n+1)\times(0,1),\; \Theta\in\left(\frac\pi 2,\frac {3\pi} 4\right)\right) >0 ~.
$$
Thus, the ergodicity of the Poisson point process leads to
\begin{equation}
\label{equa:ergo1}
\mathbf{P}\left(\exists n \ge 1,\exists x=(\xi,\Theta,V)\in\mathbf{X}\ ;\ f_{\mathbf{X}}(x)=\infty,\; \xi\in[n,n+1)\times(0,1),\; \Theta\in\left(\frac\pi 2,\frac {3\pi} 4\right)\right) =1 ~.
\end{equation}
Besides, a similar argument gives
\begin{equation}
\label{equa:ergo2}
\mathbf{P}\left(\exists n \ge 1,\exists x=(\xi,\Theta,V)\in\mathbf{X}\ ;\ f_{\mathbf{X}}(x)=\infty,\; \xi\in[-n-1,-n)\times(0,1),\; \Theta\in\left(\frac {\pi} 4, \frac {\pi} 2\right)\right) = 1 ~.
\end{equation}
From (\ref{equa:ergo1}) and (\ref{equa:ergo2}), we deduce by a simple geometric argument that a.s. there exist $x\not=x'\in\mathbf{X}$ with $f_{\mathbf{X}}(x)=f_{\mathbf{X}}(x')=\infty$ and such that $\text{Grain}(x,\infty) \cap \text{Grain}(x',\infty) \not= \emptyset$. This contradicts the hardcore property in Definition \ref{def:stopcong} and completes the proof for the unilateral line segment model.

The proofs for the bilateral line segment model and the Brownian model are simpler since, in both cases, any two given grains $\text{Grain}(x,\infty)$ and $\text{Grain}(x',\infty)$, for $x\not=x'\in\mathbf{X}$, have to overlap with probability $1$. {Recall that the difference of two independent Brownian motions is still a Brownian motion which is recurrent in dimension $2$ and then eventually hits the origin.} 
\end{dem}

\section{Proof of Theorem \ref{theo:ex}}
\label{Section:existence}

\subsection{Sketch of the Proof}\label{section_SketchTh31}

{We have to prove that a $(\lambda,\delta,\mathcal{L})$-germ-grain model  $\mathbf X$ is stopped under the assumptions \eqref{equ:suffi} and \eqref{assregular}. Following Definitions \ref{def:stopcong} and \ref{def:pog}, it means that there exists an unique stopped exploration which provides the lifetime of each growing grain until it is stopped. We use a dynamical construction letting the time runs and declaring progressively if the the grains are stopped or if they are still alive.}   

\begin{figure}[!ht]
\begin{center}
\psfrag{x1}{\small{$x_{1}$}}
\psfrag{x2}{\small{$x_{2}$}}
\psfrag{x3}{\small{$x_{3}$}}
\psfrag{x4}{\small{$x_{4}$}}
\includegraphics[width=7.5cm,height=5.5cm]{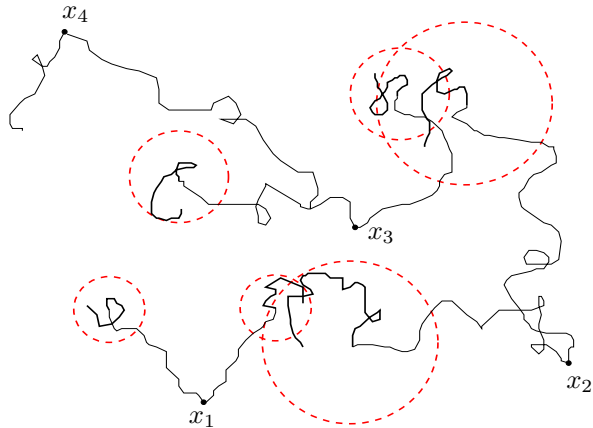}
\caption{\label{fig:obs} On this picture, the four grains (in black) associated to the marked points $x_1,x_2,x_3,x_4$, each contain two branches. At time $t$, only the grain of $x_{4}$ is already stopped (by the one associated to $x_3$). For the three other marked points still alive, we draw the red circles centered at the corresponding (six) particles with radii given by the $M_{t,t'}$'s. Remark that $x_1,x_2,x_3$ belong to the same genealogical cluster: $x_1 \sim_{t,t'} x_2$ and $x_2 \sim_{t,t'} x_3$}
\end{center}
\end{figure}

{Let us assume that the states of grains, stopped or still alive, are known at a given time $t$. During the time interval $[t,t+t']$, a particle, still alive, coming from the germ-grain $x=(\xi,k,Y)$ evolves inside a (random) disk centred at the particle with radius  given by the variable $M_{t,t'}(x)$ defined as in (\ref{MaxRadius}) by
$$
M_{t,t'}(x) := \max_{0\le l\le k-1} \, \sup_{0\le s\le t'} \, \| Y_{l}(t+s) - Y_{l}(t) \| ~.
$$
These disks are represented by red circles in Fig. \ref{fig:obs}. The main ingredient of the proof consists in establishing the absence of percolation for the random set made up with all the disks. Henceforth, the dynamic of the germ-grain model could be rigorously constructed on the time interval $[t,t+t']$ by treating independently each (finite) cluster and the corresponding grains. Indeed the dynamic for a finite number of grains is simply defined using a finite number of stopped times and a recursive procedure. Only one grain is stopped first by the others (just by checking a finite number of stopped time). Afterwards only one new grain is stopped first by the others (again by checking a finite number of stopped time) and so on... Updating the states of all grains in all finite clusters, the dynamic is now defined until the time $t+t'$ and we start a new step of the algorithm from the time $t+t'$.}

{This algorithm allows us to define the dynamic for all times if we can choose $t'>0$ small enough (uniformly with respect to $t\ge 0$) to ensure the absence of percolation for the random set of disks mentioned above. This result is proved rigorously in Proposition \ref{prop:clusterinfini}, Section \ref{sect:ProofClusterInfini}. Mainly, from Assumption \eqref{equ:suffi}, for $t'>0$ small enough (uniformly with respect to $t>0$), the mean volume of one disk given by $\pi\mathbf{E}(M_{t,t'}^2)$ is as small as we want. So, for $t'>0$ small enough, the volume fraction of space covered by all disks is so small that the percolation can not occur. The power four, at the place of the expected power two in Assumption  \eqref{equ:suffi}, is due to the multi-type branching process we use to control the percolation. Details are provided in the next Section.}

\subsection{A percolation argument}

Let us consider a germ-grain model $\mathbf{X}$ satisfying the assumptions of Theorem \ref{theo:ex}. For any times $t,t'\geq 0$ and any marked points $x=(\xi,k,Y)$, $x'=(\xi',k',Y')$ of $\mathbf{X}$ still alive at time $t$, we will write $x\sim_{t,t'} x'$ if and only if
$$
\left( \bigcup_{l=0}^{k-1} B \Big( \xi+Y_{l}(t) , \, M_{t,t'}(x) \Big) \right) \bigcap \left( \bigcup_{l=0}^{k'-1} B \Big( \xi'+Y'_{l}(t) , \, M_{t,t'}(x') \Big) \right) \, \not= \, \emptyset
$$
where $M_{t,t'}(x)$ is the random radius associated to the marked point $x=(\xi,k,Y)$ between times $t$ and $t+t'$, and defined in (\ref{MaxRadius}). Idem for $M_{t,t'}(x')$. {Roughly speaking, $x\sim_{t,t'} x'$ means that the branches from germs $x$ and $x'$ could interact during the time interval $[t,t+t']$.}

Hence, let us consider the \textbf{genealogical graph} $\mathcal{G}_{t,t'}(\mathbf{X})$ as the (non oriented) random graph whose vertex set is given by the marked points of $\mathbf{X}$ still alive at time $t$, and whose edge set is defined by the (non oriented) relation $\sim_{t,t'}$. The connected components of $\mathcal{G}_{t,t'}(\mathbf{X})$ are called \textbf{genealogical clusters}. For instance, in Fig. \ref{fig:obs}, the three marked points $x_1,x_2,x_3$ still alive at time $t$ belong to the same genealogical cluster.

Let us consider a marked point $x\in\mathbf{X}$. The genealogical cluster of $x$ in $\mathcal{G}_{t,t'}(\mathbf{X})$, denoted by $\mathcal{C}_{t,t'}(x)$, corresponds to the set of marked points still alive at time $t$ on which the evolution of the grain associated to $x$ during the time interval $[t,t+t']$ depends. The next result which is the key ingredient for the proof of Theorem \ref{theo:ex} asserts that the genealogical cluster $\mathcal{C}_{t,t'}(x)$ of $x$ is a.s. finite. Its proof is postponed to Section \ref{sect:ProofClusterInfini}.

\begin{prop}
\label{prop:clusterinfini}
There exists a (small) $t'>0$ such that for any time $t\ge 0$, the genealogical graph $\mathcal{G}_{t,t'}(\mathbf{X})$ a.s. admits only finite clusters.
\end{prop}


From now on, we set $t'$ as the time given by Proposition \ref{prop:clusterinfini}. In order to determine the evolution of the grain associated to $x$ during the time interval $[t,t+t']$, we also have to take into account the trajectories of all the grains produced before $t$ and intersecting the random set
$$
\textrm{C}_{t,t'}(x) := \bigcup_{(\xi',k',Y')\in\mathcal{C}_{t,t'}(x)} \bigcup_{l=0}^{k'-1} B \big( \xi'+Y'_{l}(t) , \, M_{t,t'}(x') \big)
$$
which is bounded thanks to the choice of $t'$. Fortunately,

\begin{lem}
\label{lem:BallFinite}
Let $t'$ given by Proposition \ref{prop:clusterinfini}. For any time $t$ and any marked point $x\in\mathbf{X}$ still alive at time $t$, the set of $y\in\mathbf{X}$ whose grain $\text{Grain}(y,t)$ until time $t$ hits $\textrm{C}_{t,t'}(x)$ is a.s. finite.
\end{lem}

In the case where the grain of $x$ would be stopped, the next result based on Assumption (\ref{assregular}) allows to avoid some pathological situations and to identify with no ambiguity its stopping grain. Proofs of Lemmas \ref{lem:BallFinite} and \ref{lem:stop} are given at the end of this section.

\begin{lem}
\label{lem:stop}
\begin{itemize}
\item[$(i)$] Two different grains can not meet simutaneously: a.s. for any $x\not= y\in\mathbf{X}$ and for any time $t$, $\text{H}(x,t)\cap\text{H}(y,t) = \emptyset$.
\item[$(ii)$] Two branches belonging to the same grain can not hit simultaneously two different grains: a.s. for any three different marked points $x,y,z$ of $\mathbf{X}$ and any $t\geq 0$,
$$
\text{H}(x,t) \cap \text{Grain}(y,t) \not = \emptyset \; \Rightarrow \; \text{H}(x,t) \cap \text{Grain}(z,t) = \emptyset ~.
$$
\end{itemize}
\end{lem}

Now, let us explain how to combine the three previous results to check that a.s. there exists a unique stopped exploration $f_{\mathbf{X}}$.

Let us start with $t=0$ and the genealogical graph $\mathcal{G}_{0,t'}(\mathbf{X})$. We treat all the genealogical clusters in the same way (and independently). Let $\mathcal{C}$ be one of them: it is finite by Proposition \ref{prop:clusterinfini}. For this first step, Lemma \ref{lem:BallFinite} is not needed. For any $x\not= y\in\mathcal{C}$, we determine by continuity of trajectories a first hitting time between the grains of $x$ and $y$ if such (first) meeting occurs. In this case, Lemma \ref{lem:stop} $(i)$ allows to know the first arrived at the meeting point among $x$ and $y$. Then, we obtain for the cluster $\mathcal{C}$ a finite sequence of hitting times: all those concerning a given grain are different by Lemma \ref{lem:stop} $(ii)$. Then a well-defined algorithm (see Section 3 of \cite{daley2014two}) identifies among these hitting times the real stopping times. Hence it determines if $x$ is stopped or not during the time interval $[0,t']$ and, if it is, by whom.

\begin{figure}[!ht]
\begin{center}
\includegraphics[width=7cm,height=6cm]{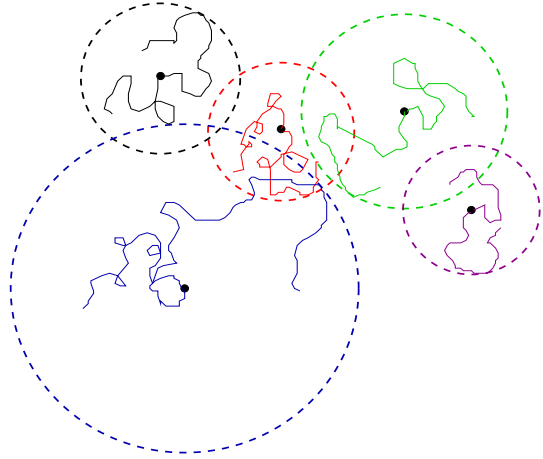}
\caption{\label{fig:bool} Here is an illustration of a genealogical cluster of $\mathcal{G}_{0,t'}(\mathbf{X})$ with $5$ marked points (with different colors) having each exactly two branches which are pictured till time $t'$. The black, green and purple marked points will be still alive at time $t'$. Only the blue and red grains meet during the time interval $[0,t']$: say the blue one stops the red one. The blue marked point will then survive until time $t'$. The parts of both red branches created after the ``red'' stopping time have to be deleted.}
\end{center}
\end{figure}

At the end of this first step, let us set $f_{\mathbf{X}}^{(1)}(x)=t'$ if $x$ is still alive at time $t'$. Otherwise $f_{\mathbf{X}}^{(1)}(x)$ is defined as its stopping time. We then obtain the set of grains already explored until time $t'$:
$$
\mathbf{G}_{1} := \bigcup_{x\in\mathbf{X}} \text{Grain}(x,f_{\mathbf{X}}^{(1)}(x)) ~.
$$

Let us continue the proof by induction. Let $k$ be a positive integer and let us assume explored and known all the grains until time $kt'$: i.e. the values $f_{\mathbf{X}}^{(k)}(x)$ for any $x$ as well as the set
$$
\mathbf{G}_{k} := \bigcup_{x\in\mathbf{X}} \text{Grain}(x,f_{\mathbf{X}}^{(k)}(x)) ~.
$$
The marked points still alive at time $kt'$ are those such that $f_{\mathbf{X}}^{(k)}(\cdot)=kt'$. For the others (those stopped before $kt'$), we can immediatly set $f_{\mathbf{X}}^{(k+1)}(x)=f_{\mathbf{X}}^{(k)}(x)$. Thus we proceed as in the first step but taking into account this time the grains produced before time $kt'$ (and using Lemma \ref{lem:BallFinite}). A similar algorithm allows to determine the evolution of grains still alive at $kt'$ until time $(k+1)t'$. If $x$ is stopped during the time interval $[kt',(k+1)t']$, then $f_{\mathbf{X}}^{(k+1)}(x)$ is defined as its stopping time. Otherwise $f_{\mathbf{X}}^{(k+1)}(x)=(k+1)t'$. Thus we update the set of explored grains:
$$
\mathbf{G}_{k+1} := \bigcup_{x\in\mathbf{X}} \text{Grain}(x,f_{\mathbf{X}}^{(k+1)}(x))
$$
which contains $\mathbf{G}_{k}$.

We can then define a function $f_{\mathbf{X}}$ valued in $(0,\infty]$ as the pointwise limit of the non-decreasing sequence $(f_{\mathbf{X}}^{(k)})_{k>0}$: for any $x\in\mathbf{X}$, $f_{\mathbf{X}}(x):=\lim \nearrow f_{\mathbf{X}}^{(k)}(x)$. In particular, $f_{\mathbf{X}}(x)$ is infinite if and only if the grain associated to $x$ is never stopped. Hence we get a unique exploration function $f_{\mathbf{X}}$ of $\mathbf{X}$ which is stopped by construction (see Definition \ref{def:stopcong}). Indeed, any two different grains do not overlap. Moreover, a stopped marked point $x$, i.e. such that $f_{\mathbf{X}}(x)<\infty$, admits only one stopping marked point by Lemma \ref{lem:stop} $(ii)$.\\

\begin{dem}{\textbf{(of Lemma \ref{lem:BallFinite})}.}
Let $t'$ given by Proposition \ref{prop:clusterinfini}. For any $t\geq 0$ and $x\in\mathbf{X}$ still alive at time $t$, the set $\textrm{C}_{t,t'}(x)$ is a.s. bounded. So, it suffices to prove that for any $R>0$, a.s. the set of $y\in\mathbf{X}$ such that $\text{Grain}(y,t)\cap B(0,R)\not=\emptyset$ is finite. We are going to prove that
$$
I := \mathbf{E} \left(\#\{y\in\mathbf{X}\ ;\ \text{Grain}(y,t)\cap B(0,R) \not= \emptyset \} \right) < +\infty ~.
$$
This holds by standard computations based on the Campbell-Mecke formula
\begin{eqnarray*}
I &\le & \lambda \left(\sum_{k=1}^{+\infty} \mathbf{P}( K = k )k\right)\int_{[0,t]}\int_{\mathbf{R}
^{2}}\int_{\mathscr{C}(\mathbf{R}_{+},\mathbf{R}^{2})}\1_{\{ \xi+U(s)\in B(0,R) \}}\Leb(ds)\Leb(d\xi)\chi(dU)\\
& \le & \lambda \mathbf{E}\left( K \right)\int_{[0,t]}\left(
\int_{\mathscr{C}(\mathbf{R}_{+},\mathbf{R}^{2})}\Leb(B(-U(s),R))\chi(dU)\right)\Leb(ds)\\
& \le & \lambda \pi R^{2}\mathbf{E}(K)t.
\end{eqnarray*}
\end{dem}

\begin{dem}{\textbf{(of Lemma \ref{lem:stop})}.}
Let us focus on Item $(i)$. Without loss of generality, we can assume that each grain contains exactly one branch. The extremity $\text{H}(.,t)$ is identified to the single point that it contains. Let us show that
$$
J := \mathbf{E}\left(\#\left\{ (x,y) \in \mathbf{X}^{2}\ ;\ x\not= y \; \mbox{ and } \; \exists t\ge 0,\ \text{H}(x,t)=\text{H}(y,t)\right\}\right) = 0 ~.
$$
By the  Campbell Mecke formula, we obtain:
\begin{eqnarray*}
J & = & \lambda^{2} \int \1_{\{\exists t\ge 0\ ;\ \xi_{1}+U_{1}(t)=\xi_{2}+U_{2}(t)\}}\left(\Leb\otimes\chi\right)^{2}\left(d(\xi_{1},U_{1}),d(\xi_{2},U_{2})\right)\\
& = & \lambda^{2}\int \1_{\{ \exists t\ge 0\ ;\ U_{1}(t)-U_{2}(t)=\xi_{2}-\xi_{1}\}}\left(\Leb\otimes\chi\right)^{2}\left(d(\xi_{1},U_{1}),d(\xi_{2},U_{2})\right)
\end{eqnarray*}
where $\chi$ denotes the marginal of $\mathcal{L}$. Thus by substitution,
\begin{eqnarray*}
J & = & \lambda^{2}\int \1_{\lbrace \exists t\ge 0\ ;\ (U_{1}-U_{2})(t)=u\rbrace}\Leb(du)\chi\otimes\chi\left(dU_{1},dU_{2}\right) \\
& = & \lambda^{2}\mathbf{E}\left(\Leb\left(\{ U_{1}(t)-U_{2}(t)\ ;\ t\ge 0\}\right) \right)
\end{eqnarray*}
which is equal to $0$ by Assumption \eqref{assregular} of Theorem \ref{theo:ex}.

Let $x_{i}:=(\xi^{(i)},k^{(i)},Y^{(i)})\in\mathbf{X}$ for $i=1,2,3$. If the grain associated to $x_1$ admits at least two branches and if it is stopped by the grain of $x_2$, we denote by $\tau=\tau(x_1,x_2)$ this stopping time. Otherwise, we set $\tau=0$. To prove Item $(ii)$, it is enough to show
$$
L := \mathbf{E} \left(\# \left\lbrace (x_{1},x_{2},x_{3})\in\mathbf{X}^{3}\ ;\ \begin{array}{c}
x_{1},x_{2},x_{3} \; \mbox{ are all different and }\\
\exists t \le \tau , \; \xi^{(1)}+Y^{(1)}_{1}(\tau) = \xi^{(3)}+Y^{(3)}_{0}(t)
\end{array} \right\rbrace \right) = 0 ~.
$$
The Campbell-Mecke formula gives:
\begin{eqnarray*}
L & = & \lambda^{3}\int \1_{\{\exists t\le\tau\ ;\ \xi^{(1)}+Y_{1}^{(1)}(\tau)=\xi^{(3)}+Y_{0}^{(3)}(u)\}} \left(\Leb\otimes\delta\otimes\mathcal{L}\right)^{3}
(dx_{1},dx_{2},dx_{3})\\
& = & \lambda^{3}\int_{\left(\mathbf{N}^{*}\times\mathcal{C}(\mathbf{R}_{+},\mathbf{R}^{2})\right)^{2}}\int_{(\mathbf{R}^{2})^{3}}\int_{\mathcal{C}(\mathbf{R}_{+},\mathbf{R}^{2})}\1_{\{ \exists t\le\tau\ ;\ Y_{0}^{(3)}(t)=\xi^{(1)}-\xi^{(3)}+Y_{1}^{(1)}(\tau)\}}\\
& & \chi(dY_{0}^{(3)})d\xi^{(1)}d\xi^{(2)}d\xi^{(3)})(\delta\otimes\mathcal{L})^{2}(d(k^{(1)},Y^{(1)}),d(k^{(2)},Y^{(2)})) ~.
\end{eqnarray*}
Now, assume that $Y^{(1)}, k^{(1)}$ and $Y^{(2)}, k^{(2)}$ are fixed. Let us specify that $\tau$ does not depend only on $Y^{(1)}, k^{(1)}$ and $Y^{(2)}, k^{(2)}$ but also on $\xi^{(1)} $ and $\xi^{(2)} $. Hence,
\begin{eqnarray*}
& & \int_{(\mathbf{R}^{2})^{3}}\left(\int_{\mathcal{C}(\mathbf{R}_{+},\mathbf{R}^{2})}\1_{\{ \exists t\le\tau\ ;\ Y_{0}^{(3)}(t)=\xi^{(1]}-\xi^{(3)}+Y_{1}^{(1)}(\tau)\}}\chi(dY_{0}^{(3)})\right)d\xi^{(1)}d\xi^{(2)}d\xi^{(3)} \\
 & = & \int_{\mathcal{C}(\mathbf{R}_{+},\mathbf{R}^{2})}\left(\int_{(\mathbf{R}^{2})^{2}}\left(
\int_{\mathbf{R}^{2}}\1_{\{\exists t\ge 0\ ;\ Y_{0}^{(3)}(t)=\xi^{(1]}-\xi^{(3)}+Y_{1}^{(1)}(\tau)\}}d\xi^{(3)}\right)d\xi^{(1)}d\xi^{(2)}\right)\chi(dY_{0}^{(3)})\\
& \le & \int_{\mathcal{C}(\mathbf{R}_{+},\mathbf{R}^{2})}\left(\int_{(\mathbf{R}^{2})^{2}}\Leb\left(\{ \eta\in\mathbf{R}^{2}\ ;\ \exists t\ge 0\ ;\ Y_{0}^{(3)}(t)=\eta\}\right)d\xi^{(1)}d\xi^{(2)}\right)\chi(dY_{0}^{(3)})\\
& \le & \int_{(\mathbf{R}^{2})^{2}}\left(\int_{\mathcal{C}(\mathbf{R}_{+},\mathbf{R}^{2})}\Leb\left(\{ \eta\in\mathbf{R}^{2}\ ;\ \exists t\ge 0\ ;\ Y_{0}^{(3)}(t)=\eta\}\right)\chi(dY_{0}^{(3)})\right)d\xi^{(1)}d\xi^{(2)}\\
& \le & \int_{(\mathbf{R}^{2})^{2}} \mathbf{E}\left(\Leb\left(\{ Y^{(3)}(t)\ ;\ t\ge 0\}\right) \right)d\xi^{(1)}d\xi^{(2)}
\end{eqnarray*}
which is null by Assumption \eqref{assregular} in Theorem \ref{theo:ex}. This implies $L=0$.
\end{dem}


\subsection{Proof of Proposition \ref{prop:clusterinfini}}
\label{sect:ProofClusterInfini}

Our proof is inspired by Theorem 2 in \cite{hall1985continuum}. Given a $(\lambda,\delta,\mathcal{L})$-germ grain model $\mathbf{X}$ and times $t\ge 0$, $t'>0$, we define a Poisson point process $\mathbf{Y}$ on $\mathbf{R}^{2}\times\mathbf{N}^{*}\times\left(\mathbf{R}^{2}\times\mathbf{R}^{*}_{+} \right)^{\mathbf{N}}$ as follows: any marked point $(\cdot,\cdot,(Y_i(t):t\geq 0)_{i\geq 0})$ is replaced with $(\cdot,\cdot,(Z_i^{(1)},Z_i^{(2)})_{i\geq 0})$ where
$$
Z_{i}^{(1)} = Y_{i}(t) \; \mbox{ and } \; Z_{i}^{(2)} = \sup_{0\le s\le t'} \|Y_{i}(t+s)-Y_{i}(t) \| ~.
$$
This transformation will simplify the notations in the following. The distribution of the sequence $(Z_i^{(1)},Z_i^{(2)})_{i\geq 0}$ is denoted by $\mathscr{W}$ and its marginal is denoted by $\Gamma$. Moreover we consider the random variable
\begin{equation}
\label{eq:rayon}
\mathbf{M}=\max_{0\le i\le K-1}Z_{i}^{(2)}
\end{equation}
in place of $M_{t,t'}$. Thus, to the Poisson process $\mathbf{Y}$, we associate the Boolean model
\begin{equation}
\label{eq:bool2s}
\text{Bool}(\mathbf{Y}) := \bigcup_{(\xi,k,Y)\in\mathbf{X}} \; \bigcup_{l=0}^{k-1} B \left(\xi+Y_{l}(t) , \max_{0\le l\le k-1} \sup_{0\le s\le t'} \|Y_{l}(t+s)-Y_{l}(t) \| \right) 
\end{equation}
and the non oriented graph $\mathcal{G}(\mathbf{Y})$ whose vertex set is given by $\mathbf{Y}$ and two vertices $(\xi,k,(Z_i^{(1)},Z_i^{(2)})_{i\geq 0})$ and $(\bar{\xi},\bar{k},(\bar{Z}_i^{(1)},\bar{Z}_i^{(2)})_{i\geq 0})$ are connected by an edge whenever
$$
\left( \bigcup_{i=0}^{k-1} B \Big( \xi+Z_{i}^{(1)} , \max_{0\le i\le k-1} Z_{i}^{(2)} \Big) \right) \bigcap \left( \bigcup_{i=0}^{\bar{k}-1} B \Big( \bar{\xi}+\bar{Z}_{i}^{(1)} , \max_{0\le i\le \bar{k}-1} \bar{Z}_{i}^{(1)} \Big) \right) \not= \emptyset ~.
$$

\begin{lem}
\label{lem:perco}
If $\lambda \mathbf{E}(\mathbf{M}^{4}) \mathbf{E}_{\delta}(K^{2})$ is smaller than $1/16\pi$. Then $\mathcal{G}(\mathbf{Y})$ a.s. does not percolate.
\end{lem}

Proposition \ref{prop:clusterinfini} is a direct consequence of Lemma \ref{lem:perco}. Indeed, by hypotheses on $\mathbf{X}$, we can choose $t'$ small enough so that $16\pi \lambda \mathbf{E}(\mathbf{M}^{4})\mathbf{E}_{\delta}(K^{2})$ is smaller than $1$, uniformly on $t$ thanks to (\ref{equ:suffi}). Thus Lemma \ref{lem:perco} applies: $\mathcal{G}(\mathbf{Y})$ and also $\mathcal{G}_{t,t'}(\mathbf{X})$ does not percolate with probability $1$.\\

In the case where $\delta=\delta_{1}$ the Dirac measure on $1$, Lemma \ref{lem:perco} is exactly Theorem 2 in \cite{hall1985continuum}. Our proof are based on the same arguments but we take account that the number of branches per germ is random. 

\begin{dem}{\textbf{(of Lemma \ref{lem:perco})}.}
We denote by $C_{0}$ the genealogical cluster of $0$ in $\mathcal{G}(\mathbf{Y}\cup\lbrace(0,K,Z\rbrace)$. We will show that $\#C_{0}$ is almost surely finite where $\#C_{0}$ denotes the number of vertices in $C_{0}$. Without loss of generality, let us replace $\mathbf{M}$ by $\lfloor\mathbf{M}\rfloor+1 $ (in other words, we consider that $\mathbf{M}$ is distributed on $\mathbf{N}^{*}$). Therefore we assume that the radii in $\text{Bool}(\mathbf{Y}) $ have positive integer values. In the following a disk is said of type $i\in\mathbf{N}^{*}$ if its radius is equal to $i$. By the same way, a given germ $\xi$ is of type $i$ if its associated disks have the type $i$. The probability $\mathbf{P}(\mathbf{M}=i)$ is denoted by $p_{i}$.

Following the strategy in \cite{hall1985continuum}), we construct a multi-type branching process such that the number of individuals dominates stochastically the number of disks in $C_{0}$. And we show that the expected number of individuals, given that $0$ is of type $i$, is bounded for any $i\ge 1$. See Athreya and Ney (\cite{athreya2004branching} page 184) for the relevant theory of multi-type branching processes. The individuals in the branching process are disks. The individuals in the $0^{\text{th}}$ generation are the $K$ disks related to the marked point $(0,K,Z)$. Given the $N^{(n)}$ individuals $\lbrace B_{l}^{(n)}\rbrace_{1\le l\le N^{(n)}} $ in the $n^{\text{th}}$ generation, we define the $(n+1)^{\text{th}}$ generation as follows. Let $\mathbf{Y}^{(n+1)}$ be a Poisson process in $\mathbf{R}^{2}\times\mathbf{N}^{*}\times\left(\mathbf{R}^{2}\times\mathbf{R}^{*}_{+}\right)^{\mathbf{N}} $ with intensity $\lambda\Leb\otimes\delta\otimes\mathscr{W}$, independent with the previous history of the process. The individuals in the $(n+1)^{\text{th}}$ generation are the disks of $\text{Bool}(\mathbf{Y}^{(n+1)}) $ which have at least one associated disk overlapping the boundary of one disk in the $n^{\text{th}}$ generation. 

Let us introduce some useful random variables related to the branching process. For two integers $n\in\mathbf{N}$ and $i\in\mathbf{N}^{*}$, we define the random variable $N_{i}^{(n)}$ as the number of disks of type $i$ in the $n^{\text{th}}$ generation, then $N^{(n)}=\sum_{i} N_{i}^{(n)}$. The expected value of $N_{i}^{(n)}$ will be denoted by $v_{i}^{(n)}$. Then, we define:
\begin{eqnarray*}
v^{(n)} &=& \left(v_{1}^{(n)},v_{2}^{(n)},\dots,v_{i}^{(n)},\dots\right),\\
\#v^{(n)} &=& \sum_{j=1}^{+\infty}v_{j}^{(n)}.
\end{eqnarray*}
In other words, $\#v^{(n)}$ is simply the expected total number of individuals in the $n^{\text{th}}$ generation. Then, we have the following relation 
\begin{equation}\label{equa:trans}
 N_{i}^{(n+1)}=\sum_{j=1}^{+\infty}\sum_{l=1}^{N_{j}^{(n)}}\mathbf{U}_{(j,i)}^{(l)} ,
\end{equation}
where
\begin{itemize}
\item For $j\in\mathbf{N}^{*}$, the $N_{j}^{(n)}$ disks of type $j$ are ranked using the lexicographic order of their centers $(B^{(n)}_{(j,l)})_{1\le l\le N_{j}^{(n)}}$.
\item For $i,j,l\in\mathbf{N}^{*}$, $$\mathbf{U}_{(j,i)}^{(l)}=\sum_{(\xi,k,Z)\in\mathbf{Y}^{(n+1)}}
k\mathbb{1}_{\{\xi\text{ is of type }i\text{ and one of the }k\text{ disks starting from }\xi\text{ overlap }\partial B^{(n)}_{(j,l)}\}}.$$ 
\end{itemize}
The probability measure of $\mathbf{U}_{(j,i)}^{(l)}$ does not depend on $N_{j}^{(n)}$. Hence, for any $1\le l\le N_{j}^{(n)}$
$$\mathbf{E}\left(\mathbf{U}_{(j,i)}^{(l)}\right)=\mathbf{E}\left(\sum_{(\xi,k,Z)\in\mathbf{Y}}
k \1_{\{\xi\text{ is of type }i\text{ and one of the }k\text{ disks starting from }\xi\text{ overlaps }\partial B(0,j)\}}\right) .$$
In the sequel $\mathbf{E}\left(\mathbf{U}_{(j,i)}^{(l)}\right)$ is denoted by $\mu_{(i,j)}$. Then,
$$\mu_{(i,j)}
\le \mathbf{E}\left(\sum_{(\xi,k,Z)\in\mathbf{Y}}k^{2} \1_{\{\xi\text{ is of type }i\text{ and the first disk starting from }\xi\text{ overlaps }\partial B(0,j)\}}\right),$$
where the first ball of a given marked point $(\xi,k,Z)$ is the one with centre $\xi+Z^{(1)}_{1}$. 

For a given marked point $y=(\xi,k,Z)$, let us define the genre of $y$ as the couple of integers $(i,k)$ such that $i$ is the radius of any disk starting from $\xi$, and $k$ is the number of disks starting from $\xi$. It follows that:
\begin{eqnarray*}
\mu_{(j,i)}
&\le &\sum_{k=1}^{+\infty} k^{2}\mathbf{E}\left(\#\left\{y=(\xi,k,Z)\in\mathbf{Y}\ \text{of genre } (i,k)\ \text{such that }B(\xi+Z^{(1)}_{0},i)\cap B(0,j)\not=\emptyset\right\}\right),\\
& \le & \sum_{k=1}^{+\infty}k^{2}\mathbf{E}\left(\#\lbrace y=(\xi,k,Z)\in\mathbf{Y}\text{ of genre }(i,k)\text{, such that }\xi+Z_{0}^{(1)}\in B(0,i+j)\rbrace\right),
\end{eqnarray*}
We have $\Leb(B(0,i+j))= 4\pi (i+j)^{2} $. We define $\mathscr{Q}:=\mathbf{E}\left(\#\lbrace y\in\mathbf{Y}\text{ of genre }(i,k)\ ;\ \xi+Z_{0}^{(1)}\in R^{i}_{j}\rbrace\right)$ and  $\Upsilon_{i}$  the law of $Z_{0}^{(1)}$ given the event $\{\mathbf{M}=i\}$. Using the Slivnyak-Mecke formula,
\begin{eqnarray*}
\mathscr{Q} &= & \lambda \delta(\{k\}) \int_{\mathbf{R}^{2}}\left(\int_{(\mathbf{R}^{2}\times\mathbf{R}^{*}_{+})^{\mathbf{N}}}
\1_{\{ (\xi,k,Z)\text{ is of type }i\}}
\1_{\{ \xi+Z_{0}^{(1)}\in B(0,i+j)\}}\Leb(d\xi)\right)\mathscr{W}(dZ)\\
 &=& \lambda \delta (\{k\})p_{i}\int_{\mathbf{R}^{2}}\left(\int_{\mathbf{R}^{2}}
\1_{\{ \xi\in\tau_{-Z_{0}^{(1)}}(B(0,i+j))\}}\Leb(d\xi)\right)\Upsilon_{i}(dZ_{0}^{(1)}),\\
&=& \lambda \delta (\{k\})p_{i}\int_{\mathbf{R}^{2}}\Leb(B(0,i+j))\Upsilon_{i}(dZ_{0}^{(1)}),\\
&=& \lambda \delta(\{k\})p_{i}\Leb(B(0,i+j)),\\
&=& 4\pi \lambda\delta(\{k\})p_{i}(i+j)^{2}.
\end{eqnarray*}
Then
\begin{eqnarray}\label{eq:narray}
\mu_{(j,i)}  &\le & \sum_{k=1}^{+\infty} 4 k^{2}\pi\delta(\{k\})\lambda p_{i}(i+j)^{2} \nonumber\\
& \le & 4\mathbf{E}(K^{2})\pi \lambda p_{i} (i+j)^{2}.
\end{eqnarray}
Using the independence of processes $(\mathbf{Y}^{(n)})_{n}$, Equation (\ref{equa:trans}) gives us the bound   
\begin{equation}\label{equa:esp}
v_{i}^{(n+1)}\le\sum_{j=1}^{+\infty}v_{j}^{(n)}\mu_{(j,i)}.
\end{equation}
The vectors $v^{(n)}$ and $v^{(n+1)}$ are related by the following matrix product $v^{(n+1)}\le v^{(n)}\mathbf{A}$, where $\mathbf{A}=\left(\mu_{(k,l)}\right)$ is a matrix with infinite number of rows (indexed by $k$) and columns (indexed by $l$). By iteration, we obtain that $v^{(n)}\le v^{(0)}\mathbf{A}^{n}$, for all $n\ge 0.$ Let us assume that the initial individuals have the type $i\in\mathbf{N}^{*}$ (in this case, $v^{(0)}=(0,\dots,\mathbf{E}(K),\dots)$ the row vector whose $i^{\text{th}}$ element is $\mathbf{E}(K)$ and has all others elements zero), and $\mu_{i}$ is the expected total number of individuals in all generations, given that the initial individuals were of type $i$. Then, we obtain an upper bound for $\mu_{i}$
\begin{equation}\label{eq:depart}
\mu_{i}\le\mathbf{E}(K)+\mathbf{E}(K)\sum_{n=1}^{+\infty}\ \sum_{j=1}^{+\infty} a_{(i,j)}^{(n)},
\end{equation}
where $a_{(i,j)}^{(n)}$ is the $(i,j)^{\text{th}}$ coefficient of $\mathbf{A}^{n}$. The branching process is defined such that
\begin{equation}\label{eq:dominationB}
\mathbf{E}\left(\#C_{0}\ |\ \text{0 is of type i}\right)\le \mu_{i}.
\end{equation}
This general point of view is described in \cite{meester1996continuum} for example.
The rest of the proof consists to obtain a good upper bound for $\mu_{i}$. Using (\ref{eq:narray}), we obtain upper bounds for $a_{(i,j)}^{(n)}$. For $i,j\in\mathbf{N}^{*}$, we have:
\begin{eqnarray*}
a_{(i,j)}^{(2)}=\sum_{l=1}^{+\infty}a_{(i,l)}^{(1)}a_{(l,j)}^{(1)}&= &\sum_{l=1}^{+\infty}
\mu_{(i,l)}\mu_{(l,j)},\\
 &\le & \left(4\mathbf{E}(K^{2})\pi \lambda\right)^{2}p_{j}\sum_{l=1}^{+\infty}p_{l}(j+l)^{2}(i+l)^{2},\\
 & \le & \left(4\mathbf{E}(K^{2})\pi \lambda\right)^{2}p_{j}\sum_{l=1}^{+\infty}16 p_{l}(ijl)^{4},\\
 &\le& \left(16\mathbf{E}(K^{2})\pi\mathbf{E}(\mathbf{M}^{4}) \lambda\right)^{2}p_{j}j^{4}i^{4}.
\end{eqnarray*}
If, for all $i,j$, we have $a_{(i,j)}^{(n-1)}\le \left(16\mathbf{E}(K^{2})\pi\mathbf{E}(\mathbf{M}^{4})\lambda\right)^{n-1}p_{j}j^{4}i^{4}$ , it follows easily that $a_{(i,j)}^{(n)}\le \left(16\mathbf{E}(K^{2})\pi\mathbf{E}(\mathbf{M}^{4})\lambda\right)^{n}p_{j}j^{4}i^{4}$, and so the latter formula must be true for all $i,j$ and $n$, using mathematical induction. Substituting this estimate into (\ref{eq:depart}) we see that
$$\mu_{i}\le \mathbf{E}(K)+ i^{4}\mathbf{E}(K)\mathbf{E}(\mathbf{M}^{4})\sum_{n=1}^{+\infty}\left(16\mathbf{E}(K^{2})\pi\mathbf{E}(\mathbf{M}^{4}) \lambda\right)^{n}.$$
which is finite if $16\mathbf{E}(K^{2})\pi\mathbf{E}(\mathbf{M}^{4}) \lambda<1$. In this case, the expected total number of individuals in all generations (given that the initial individuals were of type $i$) is finite, and consequently $\mathbf{E}\left(\#C_{0}\ |\ \text{0 is of type i}\right)<+\infty$.
\end{dem}

\section{Proof of Theorem \ref{theo:perco}}
\label{Section_perco}

This section aims to prove that the unilateral line segment model $\mathbf{X}$ defined in Section \ref{sect:ULS} does not percolate under the following moment condition on the velocity $\mathbf{V}$:
\begin{equation}
\label{hypo:Velocity}
\exists s > 1 , \, \mathbf{E} \Big( \exp(\mathbf{V}^{s}) \Big) < +\infty ~.
\end{equation}

\subsection{Sketch of the proof}

{The hypothesis (\ref{hypo:Velocity}) ensures that any line segment will be eventually stopped by another line segment. This allows us to interpret the line segment model $\mathbf{X}$ as an oriented outdegree-one graph $\mathcal{G}(\mathbf{X})$. Hence, from each marked point $x\in\mathbf{X}$ starts a unique oriented path in $\mathcal{G}(\mathbf{X})$, denoted by $\text{For}(\mathbf{X},x)$ and called the \textit{Forward set} of $x$. By taking advantage of the stationarity in distribution of the model, we can reduce the absence of percolation to the fact that a.s. any Forward set $\text{For}(\mathbf{X},x)$ is finite, which means that each $\text{For}(\mathbf{X},x)$ ends with a loop. See Section \ref{sect:Step1}.}

{Let us proceed by contradiction with assuming that from a typical marked point $\gamma$, located at the origin, may start an infinite forward path (i.e. with no loop):
\begin{equation}
\label{eq:dep}
\mathbf{P} \left( \#\text{For}(\mathbf{X}_{\gamma},\gamma) = \infty \right) > 0
\end{equation}
where $\mathbf{X}_{\gamma}:=\mathbf{X}\cup\{\gamma\}$. A marked point $x\in\mathbf{X}$ is an \textit{almost looping point} if it is possible to immediately break and stop the path $\text{For}(\mathbf{X},x)$ by adding to the current configuration three points suitably chosen. See Definition \ref{DefinitionALP} and Figure \ref{fig:cicleass}.}

{Because almost looping points are all opportunities to stop the forward path at which they belong, we prove in a second step that it is impossible for a Forward set $\text{For}(\mathbf{X},x)$ to contain infinitely many almost looping points. See Section \ref{sect:Step2}.}

{On the other hand, we show that the fact of being an almost looping point is not a very demanding property. Precisely, in Section \ref{sect:Step3}, we explicit a random region $\mathcal{S}$ conducive to almost looping points (Step 3). Moreover in Section \ref{sect:Step4} we state that this random region $\mathcal{S}$ is supercritical (Step 4). Combining these two statements, we then conclude in Section \ref{sect:Step5} that the infinite typical Forward set $\text{For}(\mathbf{X}_{\gamma},\gamma)$ visits infinitely many times the region $\mathcal{S}$ and then has to contain infinitely many almost looping points, which is forbidden by Step 2.}

\subsection{Absence of forward percolation is enough (Step 1)}
\label{sect:Step1}

Recall that under (\ref{hypo:Velocity}), the unilateral line segment model $\mathbf{X}$ is a stopped germ-grain model satisfying the finite time property: see Corollary \ref{coro:models} and Proposition \ref{prop:POG}. Hence, it is natural to associate to it an oriented Poisson Outdegree-one Graph (POG) since each marked point $x\in\mathbf{X}$ points out to its unique stopping marked point. The concept of POGs is our main structural tool to deeply understand the percolation properties of stopped germ-grain models having the finite-time property.

\begin{defin}
Let $\mathbf{X}$ be the unilateral line segment model defined in Section \ref{sect:ULS}. A.s. the unique stopping marked point of any $x\in\mathbf{X}$ is denoted by $h(\mathbf{X},x)$. Hence, we associate to $\mathbf{X}$ a Poisson outdegree-one graph (POG) $\mathcal{G}(\mathbf{X})$ which is an oriented graph whose vertex set is $\mathbf{X}$ and edge set is $\{(x,h(\mathbf{X},x)) : x\in\mathbf{X}\}$. Moreover, we denote by $h_{g}(\mathbf{X},x)$ the single point contained in $\text{H}(x,f_{\mathbf{X}}(x))$. Geometrically, $h_{g}(\mathbf{X},x)$ represents the impact point of the line segment $x$ over the line segment $h(\mathbf{X},x)$.
\end{defin}

\begin{figure}[!ht]
\begin{center}
\psfrag{x}{\small{$x$}}
\psfrag{hx}{\small{$h(\varphi,x)$}}
\psfrag{y}{\small{$h(\varphi,y)$}}
\psfrag{hy}{\small{$y$}}
\psfrag{hgx}{\small{$h_{g}(\varphi,x)$}}
\psfrag{hgy}{\small{$h_{g}(\varphi,y)$}}
\begin{tabular}{cp{2cm}c}
\includegraphics[width=6cm,height=4cm]{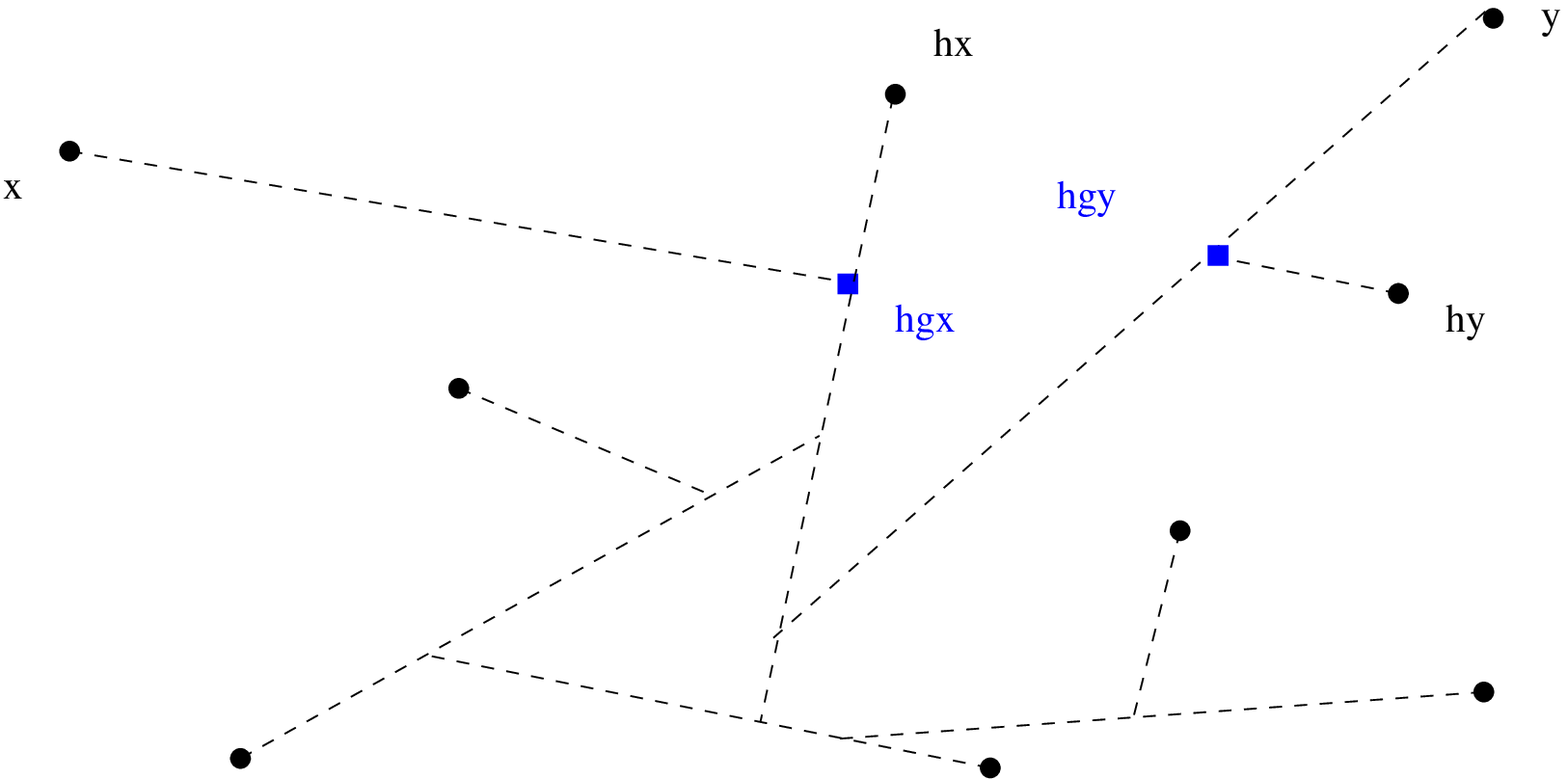} & & \includegraphics[width=6cm,height=4cm]{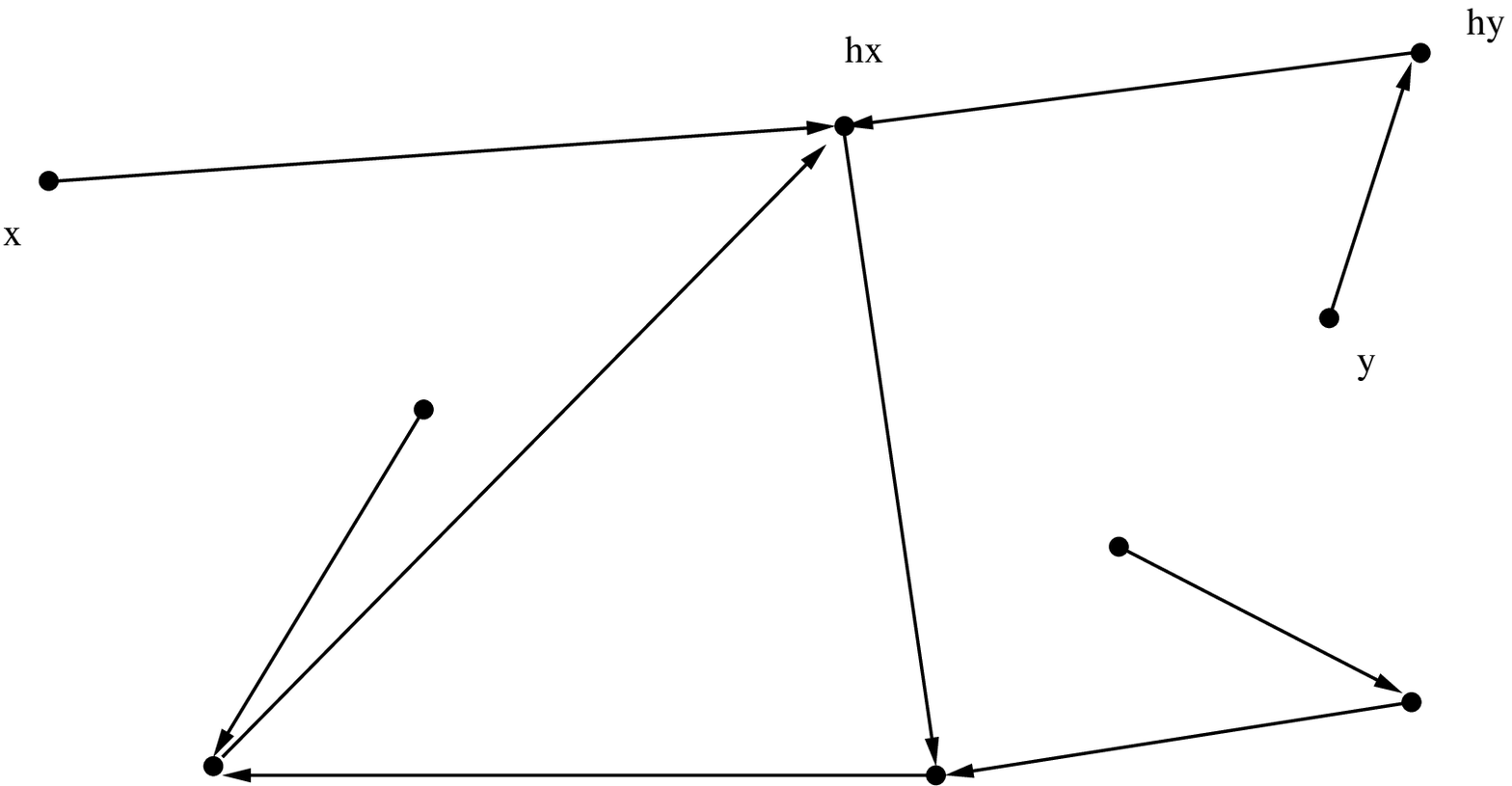}
\end{tabular}
\caption{\label{fig:segvit} Here is a finite configuration of the unilateral line segment model (to the left) viewed as a POG (to the right). The blue squares represent the impact points $h_{g}(\mathbf{X},x)$ and $h_{g}(\mathbf{X},y)$ where $x,y\in\mathbf{X}$ are marked points. On this picture, $h(\mathbf{X},x)$ is the outgoing vertex of $h(\mathbf{X},y)$ what is written as $h(\mathbf{X},x)=h(\mathbf{X},h(\mathbf{X},y))$. Remark also that the Forward sets of $x$ and $y$ contain a loop (the same one) of size $3$.}
\end{center}
\end{figure}

This formalism corresponds exactly to the concept of POGs developed in \cite{coupier2016absence} whose we recall the main notions right now. With probability $1$, for any marked point $x\in\mathbf{X}$, the \textit{Forward set} $\text{For}(\mathbf{X},x)$ of $x$ in $\mathbf{X}$ is defined as the sequence of outgoing vertices starting at $x$:
$$
\text{For}(\mathbf{X},x) := \left\{ x,h(\mathbf{X},x) , h(\mathbf{X},h(\mathbf{X},x)) , \dots \right\} ~.
$$
By construction, the forward set $\text{For}(\mathbf{X},x)$ is a branch of the POG possibly infinite. The \textit{Backward set} $\text{Back}(\mathbf{X},x)$ of $x$ in $\mathbf{X}$ is made up with all the marked points $y$ having $x$ in their Forward set, i.e.
$$
\text{Back}(\mathbf{X},x) := \{ y \in \mathbf{X} : \, x \in \text{For}(\mathbf{X},y) \} ~.
$$
The Backward set $\text{Back}(\mathbf{X},x)$ then admits a tree structure with root $x$. The Forward and Backward sets of $x$ may overlap; they (at least) contain $x$. Their union $\text{C}(\mathbf{X},x)$ forms the \textit{Cluster} of $x$ in $\mathbf{X}$. Although the Cluster $\text{C}(\mathbf{X},x)$ is a subset of the connected component of $x$ in the POG $\mathcal{G}(\mathbf{X})$, the absence of infinite clusters clearly implies the one of infinite connected components.

Also the outdegree-one property of $\mathcal{G}(\mathbf{X})$ forces each cluster to contain at most one \textit{loop}, i.e. a finite subset $\{y_{1},\ldots,y_{l}\} \subset \text{For}(\mathbf{X},x)$ such that for any $1\le i\le l$, $h(\mathbf{X},y_{i})=y_{i+1}$ (where the index $i+1$ is taken modulo $l$). The integer $l$ is called the \textit{size} of the loop. It is not difficult to observe that the Forward set $\text{For}(\mathbf{X},x)$ is finite if and only if it contains a loop. Hence, a finite cluster is made up with some finite trees merging on the same loop as illustrated in Figure \ref{fig:segvit}.  The notion of loop will be central in our study.

Let us finally recall a general percolation result for POGs which is a direct consequence of the mass transport principle. See Proposition 4.1 of \cite{coupier2016absence} for details.

\begin{prop}
\label{prop:mass}
Let $\mathbf{X}$ be a stopped germ-grain model satisfying the finite time property. Then,
$$
\mathbf{P}\left(\forall x\in\mathbf{X},\ \#\text{C}(\mathbf{X},x)<+\infty\right) = 1 \; \Longleftrightarrow	\; \mathbf{P}\left(\forall x\in\mathbf{X},\ \#\text{For}(\mathbf{X},x)<+\infty\right) = 1 ~.
$$
\end{prop}

Henceforth, our goal is to show that a.s. each Forward set in the unilateral line segment model contains a loop. The main result of \cite{coupier2016absence} (Theorem 3.1) provides two sufficient assumptions, namely the Loop assumption and the Shield assumption, ensuring the absence of forward percolation for POGs. This result applies to the line segment model with \textit{bounded velocities}: see Theorem 3.2 in \cite{coupier2016absence}. However, the Shield assumption which expresses a certain stabilizing property of the model, is no longer true in the context of \textit{unbounded velocities}: roughly speaking the reader may think about a very distant-- but very quick --marked point whose line segment would destroy any structure in a given domain. We then have to adapt the proofs of \cite{coupier2016absence} to unbounded velocities.

\subsection{A Forward set cannot contain infinitely many almost looping points (Step 2)}
\label{sect:Step2}

This part globally follow the strategy established in the proof of Theorem 3.1 in \cite{coupier2016absence}. For this reason, we skip the details.

In the sequel, we will use the notation $B(x,r)$ (instead of $B(\xi,r)$ where $x=(\xi,\cdot,\cdot)$) for the open Euclidean ball with center $\xi$ and radius $r>0$, and $\mathbf{X}_{\Gamma}$ for the elements of $\mathbf{X}$ whose first coordinate is in $\Gamma$.

\begin{defin}
Let $r<R$ be some positive real numbers and $K$ be a positive integer. A marked point $x\in\mathbf{X}$ is said a $(r,R,K)$-\textbf{looping point} of $\mathbf{X}$ if $\#\mathbf{X}_{B(x,R)}\le K$ and its Forward set $\text{For}(\mathbf{X},x)$ contains a loop $\{y_{1},\dots,y_{l}\}$, for some $l\geq 1$, whose center of mass belongs to $B(x,r)$.
\end{defin}

Let $\gamma:=(0,\Theta,V)$ be a typical marked point at the origin with direction $\Theta$ (uniformy distributed on $[0,2\pi]$) and velocity $V$. Let $\mathbf{X}_{\gamma}:=\mathbf{X}\cup\{\gamma\}$. {Now, assume that $\gamma$ is a $(r,R,K)$-looping point of $\mathbf{X}_{\gamma}$ (with loop $\{y_{1},\dots,y_{l}\}$) and consider the mass transport in which each element of $\text{Back}(\mathbf{X}_{\gamma},\gamma)$ sends a mass $1$ to the loop $\{y_{1},\dots,y_{l}\}$. Then the mass transport principle, based on the stationarity of the model, asserts that the incoming mass is equal to the outcoming mass in expectation. This implies that the loop $\{y_{1},\dots,y_{l}\}$ receives a finite (incoming) mass and thus that the expected size of $\text{Back}(\mathbf{X}_{\gamma},\gamma)$ is finite.}

\begin{prop}(Proposition 4.5 of \cite{coupier2016absence})
\label{PropBLP}
Any triplet $(r,R,K)$ satisfies
\begin{equation}
\label{BLP}
\mathbb{E} \Big( \#\text{Back}(\mathbf{X}_{\gamma},\gamma) \1_{ \{\gamma \text{ is a $(r,R,K)$-looping point of } \mathbf{X}_{\gamma}\}} \Big) < \infty ~.
\end{equation}
\end{prop}

Let $v\in\mathbf{R}^{2}$. The translation operator $\tau_v$ acts on $\mathbf{R}^{2}$ and $\mathbf{R}^{2}\times [0,2\pi]\times\mathbf{R}^{*}_{+}$ as follows: for any $w\in\mathbf{R}^{2}$, $x=(\xi,\Theta,V)\in \mathbf{R}^{2}\times [0,2\pi]\times\mathbf{R}^{*}_{+}$ and $A\subset\left(\mathbf{R}^{2}\times [0,2\pi]\times\mathbf{R}^{*}_{+}\right)$, we set $\tau_v(w)=v+w$, $\tau_v(x)=(\xi+v,\Theta,V)$ and $\tau_{v}(A)=\bigcup_{x\in A}\tau_{v}(x)$. An \textit{almost looping point} whose definition below is directly inspired by Definition 4.2 of \cite{coupier2016absence}, is set to become a looping point by adding some suitable marked points. In other words, the Forward set of an almost looping point can be stopped (by a loop) provided that the configuration is augmented with a finite number of well chosen points.   

\begin{defin}
\label{DefinitionALP}
Let us consider real numbers $0<r<R$, a positive integer $K\in\mathbf{N}^{*}$, a maximal velocity $W>0$ and an open ball $A\subset \left( B(0,r)\times[0,2\pi]\times[0,W] \right)^{3}$. A marked point $x\in\mathbf{X}$ is said a $(r,R,W,K,A)$-\textbf{almost looping point} of $\mathbf{X}$ if $\#\mathbf{X}_{B(x,R)}\le K$ and for any triplet $(x_{1},x_{2},x_{3})\in A_{x}$, we have:
\begin{itemize}
\item[$(i)$] $\text{For}(\mathbf{X}\cup\lbrace x_{1},x_{2},x_{3}\rbrace,x) = \lbrace x,x_{1},x_{2},x_{3}\rbrace$;
\item[$(ii)$] $\#\text{Back}(\mathbf{X}\cup\lbrace x_{1},x_{2},x_{3}\rbrace,x) \ge \#\text{Back}(\mathbf{X},x)$;
\end{itemize}
where $A_{x}=\tau_{\xi}(A)$ with $x=(\xi,\cdot)$.
\end{defin}

\begin{figure}[!ht]
\begin{center}
\psfrag{x}{\small{\hspace*{1cm}$x$}}
\psfrag{y}{\small{$h(\mathbf{X},x)$}}
\psfrag{z1}{}
\psfrag{z2}{}
\includegraphics[width=7.3cm,height=4.7cm]{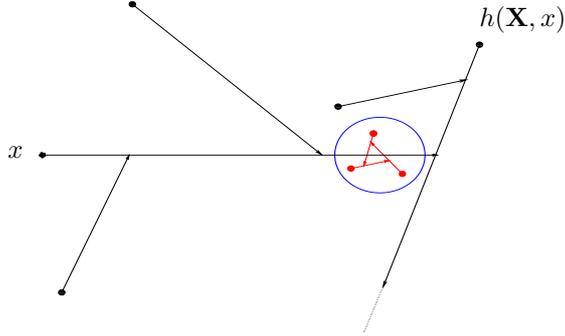}
\end{center}
\caption{\label{fig:cicleass} {On this picture, the marked point $x\in\mathbf{X}$ is an almost looping point: adding the three red points allows to stop the Forward set of $x$ without reducing its Backward set.}}
\end{figure}

Definition \ref{DefinitionALP} says that a $(r,R,W,K,A)$-almost looping point $x$ of $\mathbf{X}$ becomes a $(r,R,K+3)$-looping point of the augmented configuration $\mathbf{X}\cup\{x_{1},x_{2},x_{3}\}$ whenever the marked points $x_{1},x_{2},x_{3}$ are added in $A_{x}$; these three marked points create a loop of size $3$ which stops the line segment of $x$ as in Figure \ref{fig:cicleass}. Hence, $A_{x}$ can be understood as a suitable region to break the Forward set of $x$.

Here is the main result of this section: a typical Forward set cannot contain infinitely many almost looping points.

\begin{prop}
For any $0<r<R$, $K\in\mathbf{N}^{*}$, $W>0$ and any open ball $A$ included in $(B(0,r)\times[0,2\pi]\times[0,W])^{3}$, with probability $1$,
\begin{equation}
\label{NoInfinitelyManyALP}
\#\{ y \in \text{For}(\mathbf{X}_{\gamma},\gamma) ;\, y\text{ is a $(r,R,W,K,A)$-almost looping point of } \mathbf{X}_{\gamma} \} < \infty ~.
\end{equation}
\end{prop}

{Let us briefly recall its proof. Exploiting the outdegree-one graph structure of the model, Proposition 4.3 of \cite{coupier2016absence} says that the mean size of the Backward set of a typical almost looping point is infinite, i.e.
\begin{equation}
\label{InfiniteBackALP}
\mathbb{E} \left( \#\text{Back}(\mathbf{X}_{\gamma},\gamma) \1_{\{\gamma\text{ is a $(r,R,W,K,A)$-almost looping point of } \mathbf{X}_{\gamma} \}} \right) = \infty ~,
\end{equation}
as soon as the Forward set of a typical marked point contains an infinite number of almost looping points with positive probability. Thus, using Item $(ii)$ in Definition \ref{DefinitionALP}, we prove that (\ref{InfiniteBackALP}) implies
$$
\mathbb{E} \left( \#\text{Back}(\mathbf{X}_{\gamma},\gamma) \1_{\lbrace \gamma \text{ is a $(r,R,K+3)$-looping point for } \mathbf{X}_{\gamma}\rbrace} \right) = \infty
$$
which is forbidden by Proposition \ref{PropBLP}.}\\

\textbf{Conclusion.} {Let $A:=[0,1]^{2}\times[0,2\pi]\times\mathbf{R}^{*}_{+}$. Recall that $\lambda>0$ is the intensity of the (Poisson) point process of germs (see Section \ref{sect:ULS}) and $\gamma=(0,\Theta,V)$ is a typical marked point at the origin. Then, the stationarity of the model and the Campbell Mecke Formula give:
\begin{eqnarray*}
\mathbb{E} \left[ \sum_{x\in\mathbf{X}\cap A}  \1_{\{ \#\text{For}(\mathbf{X},x)=\infty \}} \right] & = & \int_A \mathbb{P} \big( \#\text{For}(\mathbf{X}\cup\{x\},x)=\infty \big) \, \mathbf{X} (dx) \\
& = & \lambda \, \mathbb{P} \big( \#\text{For}(\mathbf{X}_\gamma,\gamma) = \infty \big) ~.
\end{eqnarray*}}
Henceforth, proceeding by contradiction with assuming that, with positive probability, there exist infinite forward sets is equivalent to assume (\ref{eq:dep}), i.e. $\mathbf{P}(\#\text{For}(\mathbf{X}_{\gamma},\gamma)=\infty)>0$.

To get a contradiction from (\ref{eq:dep}), we are going to prove that the infinite Forward set $\text{For}(\mathbf{X}_{\gamma},\gamma)$ contains with positive probability infinitely many almost looping points. This is the goal of the next two steps.

\subsection{A conducive region $\mathcal{S}$ to almost looping points (Step 3)}
\label{sect:Step3}

Let $m$ be a positive integer tending to infinity. A marked point (or a line segment) $x=(\xi_x,\Theta_x,V_x)$ is said \textit{quick} whenever its growth velocity $V_x$ is larger than some \textit{critical velocity} $V_{c}(m)$ which will be specified later (keep in mind that $V_{c}(m)\to\infty$ with $m$). Thus, let us partition the space $\mathbf{R}^2$ into blocks $mz\oplus \Lambda_m$ of size $m$, with $z\in\mathbf{Z}^{2}$ and $\Lambda_{m}:=[-m/2,m/2)^2$. For any $z\in\mathbf{Z}^{2}$, let us define $V^{\text{max}}_{m}(z)$ as the highest velocity inside the block $mz\oplus\Lambda_{m}$:
$$
V^{\text{max}}_{m}(z) := \max \{ V_{x} : x \in \mathbf{X}_{mz\oplus\Lambda_m} \} ~.
$$
Remark that at time $1$, the extremity $\text{H}(.,1)$ of any line segment coming from the block $mz\oplus\Lambda_m$ is included in the square $mz \oplus \left[-V^{\text{max}}_{m}(z)-\frac{m}{2},V^{\text{max}}_{m}(z)+\frac{m}{2}\right]^{2}$.

Let us now define the \textit{polluted set} (by quick line segments). Recall that $\lfloor\cdot\rfloor$ denotes the floor function and $B_{\infty}(z,r) = \{ y\in\mathbf{Z}^{2} :\, \|z-y\|_{\infty} < r\}$.

\begin{defin}
\label{def:field}
Let $m,\alpha$ be some positive integers. Let us define the \textit{pollution radius} $R_{\alpha,m}(z)$ of $z\in\mathbf{Z}^{2}$ as
\begin{equation*}
R_{\alpha,m}(z) := \1_{\{ V^{\text{max}}_{m}(z) \geq V_{c}(m) \}} \left( \left\lfloor \frac{V^{\text{max}}_{m}(z)}{m} + \frac{1+\alpha}{2} \right\rfloor + 1 \right) ~.
\end{equation*}
Thus, the \textit{pollution set} $\Sigma_{\alpha,m}\subset\mathbf{Z}^{2}$ is defined as the following discrete Boolean model:
\begin{equation}
\Sigma_{\alpha,m} := \bigcup_{z\in\mathbf{Z}^{2}} B_{\infty} \left( z , R_{\alpha,m}(z) \right) ~.
\end{equation}
Finally, any vertex $z$ in $\Sigma_{\alpha,m}$ is said \textit{polluted}.
\end{defin}

It is important to note that a block $mz\oplus\Lambda_{m}$ corresponding to a polluted vertex $z\in\Sigma_{\alpha,m}$, does not necessarily contain a quick line segment but is likely to be touched by one of them before time $1$; this is the meaning of ``polluted''. Conversely, whenever $z\notin\Sigma_{\alpha,m}$, the set $mz\oplus[-\frac{\alpha m}{2},\frac{\alpha m}{2}]^{2}$ will not be touched by a quick line segment until time $1$. We will take advantage of this time interval to build shield structures outside the pollution set. Precisely, we exhibit a sequence $(\mathscr{E}_{m})_{m\ge 1}$ of local events which, outside the polluted set $\Sigma_{\alpha,m}$, acts as ``elementary shields'' and satisfies the next central property, called the ALP property. From now on, we set $\alpha=16$.

\begin{prop}{(Almost Looping Point property)}
\label{prop:ALP}
There exists a sequence $(\mathscr{E}_{m})_{m\ge 1}$ of events such that:
\begin{itemize}
\item[$(a)$] \textbf{Localization.} For any $m$, the event $\mathscr{E}_{m}$ is observable in $[-8m,8m]^{2}$ as soon as $0\notin\Sigma_{16,m}$.
\item[$(b)$] \textbf{ALP property.} Let $\mathcal{W}\subset\mathbf{Z}^{2}\!\setminus\!\Sigma_{16,m}$ such that $\mathbf{Z}^{2}\!\setminus\!\mathcal{W}$ contains at least two connected components $A_{1}$ et $A_{2}$ (w.r.t the $l_{1}$-norm) satisfying for $i$ in $\{1,2\}$,
$$
\mathcal{A}_{i} := \left( A_{i} \oplus \left[-\frac{1}{2},\frac{1}{2}\right]^{2}\right) \setminus \left( \mathcal{W} \oplus\left[-8,8\right]^{2}\right) \not= \emptyset ~.
$$
Assume that for any $z\in\mathcal{W}$, $\tau_{-mz}(\mathbf{X})\in\mathscr{E}_{m}$. Then, for any $x\in\mathbf{X}_{m\mathcal{A}_{1}}$ such that $\text{For}(\mathbf{X},x) \cap \mathbf{X}_{m\mathcal{A}_{2}} \not= \emptyset$, there exists $y\in\text{For}(\mathbf{X},x) \cap \mathbf{X}_{m\mathcal{W}\oplus\left[-8m,8m\right]^{2}}$ which is an almost looping point w.r.t. parameters $(5m,6m,V_{c}(m),K,A_y)$ for suitable $K$ and $A_y$.
\end{itemize}
\end{prop}

{Let us set
$$
\mathcal{S} = \mathcal{S}_m := \big\{ z \in \mathbf{Z}^{2} : \, z \notin \Sigma_{16,m} \; \mbox{ and } \; \tau_{-mz}(\mathbf{X}) \in \mathscr{E}_{m} \big\} ~.
$$
Thanks to the ALP property $(b)$, the set $\{mz \oplus \Lambda_{m} : z \in \mathcal{S}_m\}$ is conducive to almost looping points in the following sense: whenever it is crossed by a Forward set $\text{For}(\mathbf{X},x)$, this latter (locally) admits an almost looping point $y$.}\\

{The construction of the a sequence $(\mathscr{E}_{m})_{m\ge 1}$ and the proof of the ALP property are respectively given in Sections \ref{sect:Construction} and \ref{sect:ProofALP}. We will see in Step 4 that the probability $\mathbb{P}(\mathscr{E}_{m} \,|\, 0\notin\Sigma_{16,m})$ tends to $1$ as $m\to\infty$ and that the region $\mathcal{S}_m$ is supercritical.}

\subsubsection{Construction of $(\mathscr{E}_{m})_{m\ge 1}$}
\label{sect:Construction}

The construction of $\mathscr{E}_{m}$ follows the strategy established in Section 5 of \cite{coupier2016absence}.

Let us first introduce the triangular lattice and some related notations. Let $\Pi$ be the triangular lattice:
$$
\Pi := \left\{ a\overrightarrow{i}+b\overrightarrow{j} : \, a,b \in \mathbf{Z} \right\}
$$
where $\overrightarrow{i}:=(\sqrt{3},0)$ and $\overrightarrow{j}:=(\sqrt{3}.\cos(\frac{\pi}{3}),\sqrt{3}.\sin(\frac{\pi}{3}))$. The usual graph distance on $\Pi$ is denoted by $d_{\Pi}$. We also denote by $B^{n}(z)$ and $S^{n}(z)$ the (closed) ball and sphere with center $z$ and radius $n$ w.r.t. $d_{\Pi}$. For any $z\in\Pi$, let $\H(z)$ be the Voronoi cell of $z$ w.r.t. the vertex set $\Pi$:
$$
\H(z) := \Big\{ y\in\mathbf{R}^{2},\ \|y-z\|_{2} \le \inf_{w\in\Pi\setminus\{z\}} \|y-w\|_{2}\Big\} ~.
$$
The set $\H(z)$ is a regular hexagon centred at $z$. For any integer $n\ge 0$, let us introduce the hexagonal complex of size $n$ centred in $z$ as
$$
\H^{n}(z) = \bigcup_{y\in B^{n}(z)} \H(y) ~.
$$
For $\xi\in\mathbf{R}^{2}$, we also set $\text{Hex}^{n}(\xi)=\text{Hex}^{n}(0)+\xi$. Finally, for any integer $n\ge 1$, we define the hexagonal ring $C_{n}(0)=\H^{n}(0)\setminus \H^{n-1}(0)$ (with $\H^{0}(0)=\text{Hex}(0)$).

\begin{figure}[!ht]
\begin{center}
\psfrag{x}{\small{$x$}}
\psfrag{x'}{\small{$\text{H}_{x}$}}
\psfrag{y}{\small{$y$}}
\psfrag{B1}{\small{$B_{x}$}}
\psfrag{B2}{\small{$B_{y}$}}
\psfrag{Hex}{\small{$\text{Hex}(z)$}}
\includegraphics[width=5.7cm,height=4.6cm]{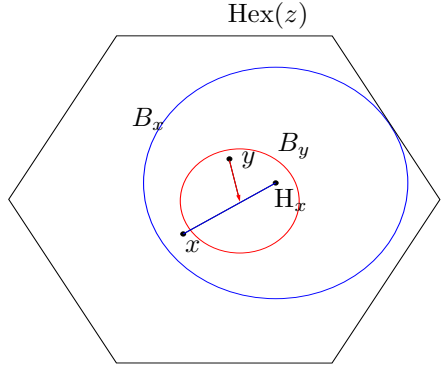}
\caption{\label{fig:local} The assumption $0\notin\Sigma_{16,m}$ allows to ``locally'' determine the line segment dynamic. Let $\text{H}_{x}:=\text{H}(x,\tau_{x})$ and $B_{x}:=B(\text{H}_x , V_{c}(m)\tau_{x})$ (with a blue circle). Let $y$ be a marked point which could stop $x$ before time $\tau_{x}$. Since $y$ is not quick, it has to belong to $\mathbf{X}_{B_{x}}$. Now let $t$ (smaller than $\tau_{x}$) be the time at which the line segment $y$ would hit the one of $x$ if it was not stopped before. Then, to check the survival of the line segment $\text{Grain}(y,\cdot)$ until $t$, we only need to observe the process $\mathbf{X}_{B_{y}}$ where $B_{y}:=B(\text{H}(y,t) , tV_{c}(m))$ (with a red circle). Then, the triangle inequality ensures that $B_{y}\subset B_{x}$. By induction, the realization of the event $\{\|\xi-h_{g}(\mathbf{X},x)\|_{2} \le V\tau_{x}\}$ only depends on $\mathbf{X}_{B_{x}}$.}
\end{center}
\end{figure}

In the sequel, we consider a positive integer $m$ and assume that $0\notin\Sigma_{16,m}$, i.e. $0$ is not polluted. Since $\text{Hex}^{4m}(0)\subset [-8 m, 8 m]^{2}$, this assumption first means that all the line segments starting from $\text{Hex}^{4m}(0)$ have a growth velocity bounded by $V_{c}(m)$. Moreover, given $z\in B^{4m}(0)$, we can also assert that no quick line segment (i.e. having a growth velocity larger than $V_{c}(m)$) may hit $\text{Hex}(z)$ before time $1$, since $\text{Hex}(z)\subset [-8 m, 8 m]^{2}$. Hence, for any $x=(\xi,\Theta,V)\in\mathbf{X}_{\text{Hex}(z)}$, there exists a time $\tau_{x}=\tau_{x}(m)\in (0,1]$ such that the realization of the event $\{\|\xi-h_{g}(\mathbf{X},x)\|_{2} \le V \tau_{x}\}$ only depends on the configuration $\mathbf{X}_{\text{Hex}(z)}$. To do it, we set $\tau_{x}$ small enough so that the Euclidean ball with center $H(x,\tau_x)$ and radius $\tau_{x} V_{c}(m)$ is included in $\text{Hex}(z)$:
$$
\tau_{x} := \sup \left\{ 0 \le t\le 1 : \, \ B(H(x,t) , tV_{c}(m)) \subset \text{Hex}(z) \right\} ~.
$$
See Figure \ref{fig:local} for an illustration.

Since then two situations may occur while observing only $\mathbf{X}_{\text{Hex}(z)}$. Either the whole exploration of the line segment of $x$ is observed until its stop before time $\tau_x$. In this case, we set $f_{m}(x):=f_{\mathbf{X}}(x)<\tau_x$. Or we can only assert that the lifetime of $x$ will be longer than $\tau_{x}$ and we set $f_{m}(x):=\tau_{x}<f_{\mathbf{X}}(x)$. In both cases, we have observed a subset $\text{Grain}(x,f_{m}(x))$ of the entire (or real) line segment $\text{Grain}(x,f_{\mathbf{X}}(x))$. Thus, we set
\begin{equation}
\label{LocalExploration}
\text{Graph}_{m}(z) := \bigcup_{x\in\mathbf{X}_{\text{Hex}(z)}} \text{Grain}(x,f_{m}(x)) ~.
\end{equation}
The crucial point is that the random set $\text{Graph}_{m}(z)$ only depends on $\mathbf{X}_{\text{Hex}(z)}$ (see also Figure \ref{fig:local}) which allows us to use later the independence property of the Poisson point process $\mathbf{X}$.\\

Let us now introduce the notion of shield hexagons.

\begin{defin}
\label{def:state1}
Let $\epsilon\in (0,1)$ and $z\in B^{4m}(0)$. The hexagon $\text{Hex}(z)$ is said \textbf{$(\epsilon,m)$-shield} for $\mathbf{X}$ if for all $a,b\in\mathbf{R}^{2}$ such that $a\notin\text{Hex}(z)$ and $b\in\epsilon\text{Hex}(z)$, we have
$$
(a,b) \cap \text{Graph}_{m}(z) \not= \emptyset ~.
$$
Moreover, for any integer $n>0$ and $\{z_{i}\}_{1\le i\le n}\subset B^{4m}(0)$, the collection $\{\text{Hex}(z_{i})\}_{1\le i\le n}$ is said \textbf{$(\epsilon,m)$-shield} for $\mathbf{X}$ if each $\text{Hex}(z_{i})$ is $(\epsilon,m)$-shield for $\mathbf{X}$.
\end{defin}

In other words, $\text{Hex}(z)$ is $(\epsilon,m)$-shield whenever its local exploration $\text{Graph}_{m}(z)$ creates a barrier (effective at time $1$ in the ring $\text{Hex}(z)\setminus \epsilon\text{Hex}(z)$). Under the assumption $0\notin\Sigma_{16,m}$, it is sufficient to observe $\mathbf{X}$ inside $\text{Hex}(z)$ to determine if $\text{Hex}(z)$ is $(\epsilon,m)$-shield for $\mathbf{X}$ or not. Hence, two hexagons $\text{Hex}(z)$ and $\text{Hex}(z')$, with $z\not=z'\in B^{4m}(0)$, are independently $(\epsilon,m)$-shield.

Furthermore, it is not difficult to convince oneself (using many small segments encircling the ring $\text{Hex}(z)\setminus \epsilon\text{Hex}(z)$, see Figure \ref{fig:pige}) that for any $\epsilon\in (0,1)$,
\begin{equation}
\label{def:pm}
p_{\epsilon,m} := \mathbf{P} \left( \text{Hex}(z) \text{ is $(\epsilon,m)$-shield } \,|\, 0 \notin \Sigma_{16,m} \right) > 0 ~.
\end{equation}
{An explicit lower bound is given in \eqref{eq:minoration} below.}

\begin{defin}
\label{def:it}
Let $\beta\in\{1,2\}$. The set $\text{Hex}^{2\beta m}(0)$ is said \textbf{$(\epsilon,m)$-shielded} for $\mathbf{X}$ if, for all $y=(\xi',\Theta',V')\in\mathbf{X}$, the segment $[\xi',h_{g}(\mathbf{X},y)]$ cannot overlap simultaneously the sets $\text{Hex}^{\beta m}(0)$ and $\mathbf{R}^{2}\setminus\text{Hex}^{2\beta m}(0)$.
\end{defin}

In other words, the set $\text{Hex}^{2\beta m}(0)$ is $(\epsilon,m)$-shielded provided no line segment crosses completely the ring $\text{Hex}^{2\beta m}(0)\setminus \text{Hex}^{\beta m}(0)$, from inside or from outside.\\

In the sequel, for $\beta\in\lbrace 1,2\rbrace$, we introduce an event $E_{m}^{(\beta)}$, whose construction is directly inspired by Section 5 in \cite{coupier2016absence}, on which $\text{Hex}^{2\beta m}(0)$ is $(\epsilon,m)$-shielded (see Lemma \ref{lem:franc}).

To do it, we need some extra notations. Let $\eta\in\partial\text{Hex}^{\beta m}(0)$ where $\partial\Lambda$ denotes the topological boundary of $\Lambda\subset\mathbf{R}^{2}$. For any $v\in[0,1]$, we define the (semi-infinite) ray starting from $\eta$ in the direction $\overrightarrow{v}:=(\cos(2\pi v),\sin(2\pi v))$ by $l(\eta,\overrightarrow{v}):=\{\eta+t\overrightarrow{v}, t\geq 0\}$. Thus, we denote by $\mathscr{L}^{m}$ the set of rays $l(\eta,\overrightarrow{v})$ coming from $\partial\text{Hex}^{\beta m}(0)$ which do not overlap the topological interior of $\text{Hex}^{\beta m}(0)$:
$$
\mathscr{L}^{m} := \left\{ l(\eta,\overrightarrow{v}) : \, (\eta,v) \in \partial\text{Hex}^{\beta m}(0) \times [0,1] \, \mbox{ and } \, l(\eta,\overrightarrow{v}) \cap \text{Int}(\text{Hex}^{\beta m}(0)) = \emptyset \right\} ~.
$$
For any ray $l\in\mathscr{L}^{m}$, the set of hexagons included in $\text{Hex}^{2\beta m}(0)\setminus\text{Hex}^{\beta m}(0)$ and crossed by $l$ is
$$
\text{Cross}(l) := \left\{ \text{Hex}(z) ,\ \beta m+1\le d_{\Pi}(0,z)\le 2\beta m \, \mbox{ and } \, l\cap\text{Hex}(z) \not= \emptyset \right\} ~.
$$
This set can be partitioned into different floors; $\text{Cross}_{i}(l)$ denotes the set of hexagons of $\text{Cross}(l)$ included in $C_{i}(0)$ for any $\beta m+1\leq i\leq 2\beta m$. Besides, let us remark that, for any ray $l\in\mathscr{L}^{m}$, there exists an index $\beta m+1\leq i(l)\leq 2\beta m$ such that, for all $i(l)\le i \le 2m$, $\text{Cross}_{i}(l)$ contains at most three hexagons (when the ray $l$ is almost parallel to the side of $\text{Hex}^{\beta m}(0)$ from which it starts, it may cross a large number of hexagons included in the same hexagonal ring $C_{i}(0)$ for $i$ close to $\beta m+1$).

The set $\text{Cross}(l)$ is said $(\epsilon,m)$\textbf{-uncrossable} for $\mathbf{X}$ whenever one can find two consecutive floors $\text{Cross}_{i}(l)$ and $\text{Cross}_{i+1}(l)$, for some index $i(l)\leq i\leq 2\beta m-1$, which are both $(\epsilon,m)$-shield for $\mathbf{X}$. Now, we can set the event $E^{\beta}_{m}(\epsilon)$ as
\begin{equation}
\label{Emn}
E_{m}^{(\beta)}(\epsilon) := \bigcap_{l\in\mathscr{L}^{m}} \left\{ \text{Cross}(l)\text{ is $(\epsilon,m)$-uncrossable for }\mathbf{X} \right\} ~.
\end{equation}

\begin{lem}
\label{lem:franc}
There exists $\epsilon\in (0;1)$ (close to 1) such that, for any $\beta\in \{1,2\}$, $\text{Hex}^{2\beta m}(0)$ is a.s. $(\epsilon,m)$-shielded for $\mathbf{X}$ on the event $E^{(\beta)}_{m}(\epsilon)\cap \{0\notin\Sigma_{16,m}\}$.
\end{lem}

\begin{dem}
Assume first that $0\notin\Sigma_{16,m}$. There is no quick line segment in $\text{Hex}^{2\beta m}(0)$ and the quick line segments from the outside $\text{Hex}^{2\beta m}(0)$ are too far to hit $\text{Hex}^{2\beta m}(0)$ before time $1$. Hence the local explorations $\text{Graph}_{m}(\cdot)$ defined in (\ref{LocalExploration}) and involved by the event $E^{(\beta)}_{m}(\epsilon)$ are realized without being disturbed by the quick line segments. Thus, the proof of Proposition 5.2 of \cite{coupier2016absence} shows that there exists $\epsilon_{\beta}$ sufficiently close to $1$ such that a.s. on the event $E^{(\beta)}_{m}(\epsilon_{\beta})$, any ray $l$ in $\mathscr{L}^{m}$ is obstructed, i.e. $l$ necessarily hits a local exploration $\text{Graph}_{m}(z)$ for some $z$ in $\text{Hex}^{2\beta m}(0)\setminus\text{Hex}^{\beta m}(0)$. Hence, the same holds for any line segments. The proof ends with $\epsilon:=\max\{\epsilon_{1},\epsilon_{2}\}$.
\end{dem}

From now on, we merely write $E^{(\beta)}_{m}$ instead of $E^{(\beta)}_{m}(\epsilon)$ and $p_{m}$ instead of $p_{\epsilon,m}$, where $\epsilon$ is given by Lemma \ref{lem:franc}. {Finally, the event $\mathscr{E}_{m}$ occurring in Proposition \ref{prop:ALP} is defined by
$$
\mathscr{E}_{m} := E_{m}^{(1)} \cap E_{m}^{(2)} ~.
$$}

\subsubsection{Proof of the ALP property}
\label{sect:ProofALP}

Let $m$ be a positive integer. Assuming that $0\notin\Sigma_{16,m}$, the construction developed above ensures that the event $\mathscr{E}_{m}$ only involves local explorations $\text{Graph}_{m}(z)$ with $z\in B^{4m}(0)$, i.e. with $\text{Hex}(z)\subset\text{Hex}^{4m}(0)$. So Item $(a)$ \textbf{Localization} of Proposition \ref{prop:ALP} follows from the inclusion $\text{Hex}^{4m}(0)\subset\left[-8m,8m\right]^{2}$. The goal of this section is to prove Item $(b)$ \textbf{ALP property} of Proposition \ref{prop:ALP}.

For this purpose, let us consider three disjoint subsets $\mathcal{W}, A_{1}, A_{2}$ of $\mathbf{Z}^{2}$ and $\mathcal{A}_{1}, \mathcal{A}_{2}$ as in Proposition \ref{prop:ALP}. We assume that $\tau_{-mz}(\mathbf{X})\in\mathscr{E}_{m}$ for any $z\in\mathcal{W}$. Thus, let us consider $x=(\xi,\Theta, V)\in\mathbf{X}_{m\mathcal{A}_{1}}$ such that $\text{For}(\mathbf{X},x) \cap \mathbf{X}_{m\mathcal{A}_{2}} \not= \emptyset$. Our goal is to identify an almost looping point $y\in\text{For}(\mathbf{X},x) \cap \mathbf{X}_{m\mathcal{W}\oplus[-8m,8m]^{2}}$ (with suitable parameters).

Let us write the Forward set $\text{For}(\mathbf{X},x)$ as a sequence $(x_{i})_{i\ge 0}$ with $x_{0}=x$ and for all index $i$, $x_{i+1}=h(\mathbf{X},x_{i})$. We also set $x_{i}:=(\xi_{i},\Theta_{i},V_{i})$ for any $i$. By assumption, there exists an index $n\ge 1$ such that $x_{n}\in \mathbf{X}_{m\mathcal{A}_{2}}$ and also an index $1\le i\le n$ and a vertex $u\in\mathcal{W}$ such that $[\xi_{i},h_{g}(\mathbf{X},x_{i})]$ overlaps $\text{Hex}^{m}(mu)$. The hypothesis $\tau_{-mu}(\mathbf{X})\in E_{m}^{(1)}$ actually implies that $[\xi_{i},h_{g}(\mathbf{X},x_{i})]$ is completely included in $\text{Hex}^{2m}(mu)$ (with Lemma \ref{lem:franc} and Definition \ref{def:it}).

\begin{figure}[!ht]
\begin{center}
\psfrag{y1}{\small{$z_{1}$}}
\psfrag{y2}{\small{$z_{2}$}}
\psfrag{y3}{\small{$z_{3}$}}
\psfrag{y4}{\small{$z_{4}$}}
\psfrag{xi}{\small{$x_{i}$}}
\includegraphics[width=7cm,height=6cm]{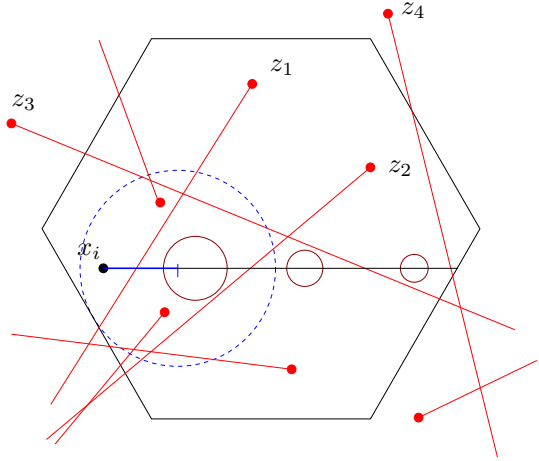}
\caption{\label{fig:last} On this picture, $l=4$ and $j=2$. The three balls $B(w_{k},r_{k})$ (for $2\le k\le 4$) are drawn in brown.}
\end{center}
\end{figure}

Let us prove that the vertex $x_{i}$ is an almost looping point of $\mathbf{X}$. Let $z_{i}$ be the element of $mz\oplus B^{2m}(0)$ such that $\xi_{i}\in\text{Hex}(z_{i})$. Recall that $\text{Grain}(x_i,f_{m}(x_i))$ is the subset that we can only observe through $\mathbf{X}_{\text{Hex}(z_{i})}$ (see Figure  \ref{fig:local}) of the true line segment $\text{Grain}(x_i,f_{\mathbf{X}}(x_i))$. Moreover, we define $\text{Grain}(x_{i},F_{m}(x_{i}))$ as the longest grain from $x_i$ remaining inside $\text{Hex}^{2m}(mu)$ (without interaction with other marked points, see Figure \ref{fig:last}). The inclusion $[\xi_{i},h_{g}(\mathbf{X},x_{i})]\subset\text{Hex}^{2m}(mu)$ forces the inequalities
$$
f_{m}(x_{i}) \leq f_{\mathbf{X}}(x_{i}) \leq F_{m}(x_{i})
$$
which means that $h_{g}(\mathbf{X},x_{i})$ lies somewhere on the segment $[\text{H}(x_{i},f_{m}(x_{i})) , \text{H}(x_{i},F_{m}(x_{i}))]$. Thus we have to identify a small ball $A_{x_i}$ corresponding to Definition \ref{DefinitionALP}, i.e. a suitable region to break the Forward set of $x_i$ without reducing its Backward set. So, this small ball $A_{x_i}$ has to be located close to (or just before) $h_{g}(\mathbf{X},x_{i})$. A difficulty appears at this stage: the hypothesis $\tau_{-mu}(\mathbf{X})\in E_{m}^{(2)}$ ensures that $x_{i+1}=h(\mathbf{X},x_{i})$ is in $\text{Hex}^{4m}(mu)$ but the observation of $\mathbf{X}$ only through $mu\oplus[-8m,8m]^{2}$ does not guarantee to identify $x_{i+1}$ and then the location of $h_{g}(\mathbf{X},x_{i})\in \text{Hex}^{2m}(mu)$. So we have to deal with a finite number of candidates for $x_{i+1}$. This is the reason why we consider
$$
\text{Ray}(u,m,x_{i}) := \bigcup_{x'\in\mathbf{X}_{\text{Hex}^{4m}(mz)}\setminus\lbrace x_{i}\rbrace}\ \text{Grain}(x',\infty) ~.
$$
A marked point $x'\in\mathbf{X}_{\text{Hex}^{4m}(mz)}\setminus\lbrace x_{i}\rbrace $ is a potential stopping marked point for $x_{i}$ if and only if $\text{Grain}(x',\infty)$ overlaps $[\text{H}(x_{i},f_{m}(x_{i})) , \text{H}(x_{i},F_{m}(x_{i}))]$. Then, there exists a random integer $l>0$ and times $0<t_{1}<\dots<t_{l}\leq F_{m}(x_{i})$ such that $\text{H}(x_{i},t_{1}),\ldots,\text{H}(x_{i},t_{l})$ are the only possible locations for $h_{g}(\mathbf{X},x_{i})$ which are created by $z_{1},\dots,z_{l}\in\mathbf{X}_{\text{Hex}^{4m}(mu)}\!\setminus\!\{x_{i}\}$ so that for any $1\leq j\leq l$, $\text{H}(x_{i},t_{j})\in\text{Grain}(z_{j},\infty)$. See Figure \ref{fig:last}.

The stopping marked point of $x_{i}$ in $\mathbf{X}$ belongs to $\lbrace z_{1},\dots,z_{l}\rbrace$. Moreover, if $t_{k}<f_{m}(x_{i})$, the marked point $z_{k}$ is not able to stop the line segment starting from $x_{i}$. Then, we set
$$
j := \min\left\{ k\in [\![1,l ]\!]\ ;\ t_{k}\ge f_{m}(x_{i})\right\} ~.
$$
The integer $j$ is well defined because we know that $h(\mathbf{X},x_{i})\in\text{Hex}^{4m}(mu)$, so there exists at least one candidate. So, $h(\mathbf{X},x_{i})\in\lbrace z_{j},\dots,z_{l}\rbrace$. To check which point is the real stopping point, we may have to control the state of $\mathbf{X}$ outside of $mu\oplus[-8m,8m]^{2}$. As we do not want to do such exploration in order to preserve the locality of the event $\mathcal E_m$, we identify a finite number of candidates for the "looping region" in Definition \ref{DefinitionALP}. 

Let us fix $k\in [\![j,l ]\!]$ and denote by $w_{k}$ the middle of the segment $\left[\text{H}(x_{i},t_{k-1}),\text{H}(x_{i},t_{k})\right]$ (with $t_{0}=0$ and $\text{H}(x_{i},0)=\xi_{i}$). There exists $r_{k}>0$ such that
\begin{equation}
\label{eq:back}
B(w_{k},r_{k})\cap\text{Ray}(u,m,x_{i})=\emptyset 
\end{equation}
and
\begin{equation}\label{eq:for}
r_{k}<\|\xi_{i}-w_{k}\|_{2} ~.
\end{equation}
In the case where $x_{i+1}$ is $z_{k}$, the obtained ball $B(w_{k},r_{k})$ is then a suitable region in which we could create an obstacle for the growing segment $x_{i}$ without reducing its Backward set. Precisely, we add a triplet of marked points in $B(w_{k},r_{k})$ which shapes a triangle of stopped line segments. This triangle will be active before the arrival of the grain from $x_{i}$ (\ref{eq:for}) and do not break any other line (\ref{eq:back}). Consider the set $\mathcal{A}^{(k)}_{x_{i}}\subset \left(B(w_{k},r_{k})\times[0,2\pi]\times[V_{i},V_{c}(m)]\right)^{3}$, where $V_i$ denotes the growth velocity of $x_i$, such that for all triplets $(a_{0},a_{1},a_{2})\in\mathcal{A}^{(k)}_{x_{i}}$:
\begin{itemize}
\item[(i)] $h(\mathbf{X}\cup\{a_{0},a_{1},a_{2}\},\  a_{l})=a_{l+1}$ for $l=0,1,2$ (where the index $l+1$ is taken modulo 3).
\item[(ii)] The triangle defined by the vertices $h_{g}(\mathbf{X}\cup\{a_{0},a_{1},a_{2}\},\  a_{l})$, $l=0,1,2$, is included in $B(w_{k},r_{k})$ and contains $w_{k}$.
\end{itemize}
It is not difficult to see that $\mathcal{A}^{(k)}_{x_{i}}$ contains a non-empty open ball $A^{(k)}_{x_{i}}\subset (B(w_{k},r_{k})\times[0,2\pi]\times[V_{i},V_{c}(m)])^{3}$. In the case where $x_{i+1}$ is $z_{k}$, by (\ref{eq:for}), (i) and (ii), any add of triplet $(a_{0},a_{1},a_{2})\in A^{(k)}_{x_{i}}$ forces the growing segment from $x_{i}$ to hit the loop created by $a_{0},a_{1},a_{2}$:
$$
\text{For}(\mathbf{X}\cup\lbrace a_{0},a_{1},a_{2}\rbrace , x_{i}) = \lbrace x_{i},a_{0},a_{1},a_{2} \rbrace ~.
$$
Moreover, condition (\ref{eq:for}) in addition with (i) and (ii) imply that no growing segment except $x_{i}$ is changing by the added marked points $\lbrace a_{0},a_{1},a_{2}\rbrace$:
$$
\text{Back}(\mathbf{X}\cup\lbrace a_{0},a_{1},a_{2} \rbrace,x_{i}) = \text{Back}(\mathbf{X},x_{i}) ~.
$$
Then, we have proved that $x_{i}$ is an almost looping point with "looping ball" $A^{(k)}_{x_{i}}$ in the case where $x_{i+1}$ is $z_{k}$. More precisely, it is not difficult to check that $x_i$ is a $(5m,6m,V_{c}(m),K,\cdot)$-almost looping point for $K$ large enough.

It is worth pointing out here that, even if the ``true looping ball'' $A_{x_{i}}$ is not precisely located, the candidates $A_{x_{i}}^{(k)}$ for $j\leq k\leq l$, only depend on the process $\mathbf{X}$ inside $\text{Hex}^{4m}(mu)\subset [-8m,8m]^{2}$. Consequently, the radius of the ball $A_{x_{i}}$ is observable w.r.t. $\mathbf{X}\cap [-8m,8m]^{2}$.

\subsection{$\mathcal{S}$ is supercritical (Step 4)}
\label{sect:Step4}

Let us set the critical velocity as follows:
$$
V_{c}(m) := \left( \log(m^{3}) \right)^{1/s} ~,
$$
where $s>1$ is given by (\ref{hypo:Velocity}). This choice comes from the following compromise. On the one hand, $V_{c}(m)$ has to tend to infinity so that the probability for a given block $mz\oplus\Lambda_{m}$ to contain at least one quick line segment tends to $0$ (Lemma \ref{lemme_quicksegments}). On the other hand, the construction of the events $(\mathscr{E}_{m})_{m\ge 1}$ requires that $V_{c}(m)$ increases very slowly (see Lemma \ref{lem:lim}). Satisfying both conditions needs a strong moment hypothesis on the speed distribution, namely (\ref{hypo:Velocity}).

\begin{lem}
\label{lemme_quicksegments}
For any vertex $z\in\mathbf{Z}^{2}$,
$$
\lim_{m\to\infty} \mathbf{P} \left( V^{\text{max}}_{m}(z) \geq V_{c}(m) \right) = 0 ~.
$$
\end{lem}

\begin{lem}
\label{lem:lim}
Given $\beta\in\{1,2\}$,
$$
\lim_{m\rightarrow +\infty} \mathbf{P} \Big( E^{(\beta)}_{m} \, | \, 0 \notin \Sigma_{16,m} \Big) = 1 ~.
$$
\end{lem}

Lemmas \ref{lemme_quicksegments} and Lemma \ref{lem:lim} (which will be proved at the end of the current section) are the main ingredients to prove that the region
$$
\mathcal{S} = \mathcal{S}_m = \big\{ z \in \mathbf{Z}^{2} : \, z \notin \Sigma_{16,m} \; \mbox{ and } \; \tau_{-mz}(\mathbf{X}) \in \mathscr{E}_{m} \big\}
$$
(conducive to almost looping points, thanks to Section \ref{sect:Step3}) is supercritical in the following sense:

\begin{prop}
\label{prop:NoPerco}
For any $m$ large enough, the set $\mathbf{Z}^{2}\!\setminus\!\mathcal{S}_m$ does not percolate with probability $1$ w.r.t the $l_{1}$-norm.
\end{prop}

\begin{dem}
Let us realize the Poisson point process $\mathbf{X}$ as the union
$$
\mathbf{X} = \mathbf{X}_{\text{quick}}^{(m)} \sqcup \mathbf{X}_{\text{slow}}^{(m)}
$$
of two independent Poisson point processes $\mathbf{X}_{\text{quick}}^{(m)}$ and $\mathbf{X}_{\text{slow}}^{(m)}$ with respective intensities $\lambda \mathbf{P}(V\geq V_{c}(m)) \otimes \Xi \otimes \mathscr{L}(V | V\geq V_{c}(m))$ and $\lambda \mathbf{P}(V<V_{c}(m)) \otimes \Xi \otimes\mathscr{L}(V | V<V_{c}(m))$. For the notations $\lambda,\Xi,V$ the reader may refer to Section \ref{sect:ULS}. Thus, let us consider the random field $\zeta:=\{\zeta^{(m)}_{z} , z\in\mathbf{Z}^{2}\}$ where
$$
\zeta^{(m)}_{z} := \1_{\{ \tau_{-mz}(\mathbf{X}^{(m)}_{\text{slow}}) \notin \mathscr{E}_{m} \}} ~.
$$
This is a stationary site percolation model with parameter $q_m:=\mathbf{P}(\mathbf{X}^{(m)}_{\text{slow}} \notin \mathscr{E}_{m})$ which is $16$-dependent, i.e. the r.v.'s $\zeta^{(m)}_{z}$ and $\zeta^{(m)}_{z'}$ are independent whenever $\|z-z'\|_{\infty}>16$. Since $\mathbf{X} \cap (mz \oplus \Lambda_{m})$ and $\mathbf{X}_{\text{slow}}^{(m)} \cap (mz \oplus \Lambda_{m})$ are equal as soon as $z$ is not polluted, we get
$$
q_m = \mathbf{P} \big( \mathbf{X}^{(m)}_{\text{slow}} \notin \mathscr{E}_{m} \big) = \mathbf{P} \big( \mathbf{X} \notin \mathscr{E}_{m} \,|\, 0 \notin \Sigma_{16,m} \big) \, \to \, 0
$$
as $m\to\infty$ by Lemma \ref{lem:lim}. Then, a classical stochastic domination result due to Liggett et al \cite{liggett1997domination} allows to stochastically dominate the (dependent) field $\zeta$ by an independent site percolation model $\xi:=\{\xi^{(m)}_{z} , z\in\mathbf{Z}^{2}\}$ with parameter $f(q_{m})$ tending to $0$ with $m$. In particular, the random set $\mathbf{Z}^{2}\!\setminus\!\mathcal{S}_m$ is included in $\Sigma_{16,m} \cup \{z : \xi^{(m)}_{z}=1\}$.

It is then sufficient to prove that this set does not percolate for $m$ large enough. To do it, it is useful to remark that $\Sigma_{16,m} \cup \{z : \xi^{(m)}_{z}=1\}$ can be viewed as the following discrete Boolean model
\begin{equation}
\label{DiscreteBoolModel}
\bigcup_{z\in\mathbf{Z}^{2}} B_{\infty} \left( z , R_{16,m}'(z) \right) ~,
\end{equation}
where $R_{16,m}'(z) := \max\big( R_{16,m}(z) , \xi^{(m)}_{z} \big)$, for any vertex $z$. Indeed, the collections $\{R_{16,m}(z) , z\in\mathbf{Z}^{2}\}$ and $\{\xi^{(m)}_{z} , z\in\mathbf{Z}^{2}\}$ are each i.i.d. families of r.v.'s and they are also independent from each other since the $R_{16,m}(z)$'s (and the polluted set $\Sigma_{16,m}$) only depends on $\mathbf{X}_{\text{quick}}^{(m)}$ whereas the fields $\zeta$ (and then $\xi$) only depends on $\mathbf{X}_{\text{slow}}^{(m)}$.

A discrete version of Theorem 2.1 of J.-B. Gou\'er\'e \cite{gouere2008subcritical} asserts that the (discrete) Boolean model defined in (\ref{DiscreteBoolModel}) is subcritical provided the mean volume of a ball is small enough. This is the reason why we are going to prove that:
\begin{equation}
\label{equa:mom}
\lim_{m\rightarrow +\infty} \mathbf{E} \Big( R'_{16,m}(0) ^{2} \Big) = 0 ~.
\end{equation}
Since $f(q_m)\to 0$, it is enough to prove that $\mathbf{E}(R_{16,m}(0)^{2})$ tends to $0$ with $m$. This immediately follows from (\ref{EspTo0-1}) and (\ref{EspTo0-2}) below. First, by definition of the radius $R_{16,m}(0)$ and using the Poisson distribution of $\mathbf{X}$, we can write:
\begin{eqnarray}
\label{EspTo0-1}
\mathbf{P} \big( R_{16,m}(0) = n \big) & \leq & \mathbf{P} \big( V_{m}^{\text{max}}(0) \geq m (n-19/2) \big) \nonumber \\
& \leq & 1 - \exp \Big( -\lambda m^{2} \mathbf{P} \big(\mathbf{V} \ge m (n-19/2) \big) \Big) \nonumber \\
& \leq & \frac{\lambda \mathbf{E}(\mathbf{V}^4)}{m^2 (n-\frac{19}{2})^4} ~,
\end{eqnarray}
for any $n\geq 10$. Thus, for the small values of $n$, we use Lemma \ref{lemme_quicksegments}:
\begin{equation}
\label{EspTo0-2}
\mathbf{P} \big( R_{16,m}(0) = n \big) \leq \mathbf{P} \left( V^{\text{max}}_{m}(z) \geq V_{c}(m) \right) \, \to \, 0 ~.
\end{equation}
\end{dem}

\begin{dem}{ (Proof of Lemma \ref{lemme_quicksegments})}
It is enough to write, using stationarity of the model and the Markov inequality:
\begin{eqnarray*}
\mathbf{P} \left( V^{\text{max}}_{m}(z) \geq V_{c}(m) \right) & = & 1 - \exp \left( -\lambda m^{2} \mathbf{P}( \mathbf{V} \geq V_{c}(m)) \right) \\
& \leq & \lambda m^{2} \mathbf{P}( \mathbf{V} \geq V_{c}(m)) \\
& \leq & \lambda \frac{\mathbf{E}(\exp(\mathbf{V}^{s}))}{m} ~.
\end{eqnarray*}
\end{dem}

\begin{figure}[!ht]
\begin{center}
\psfrag{x1}{\small{$x_{1}$}}
\psfrag{x2}{\small{$x_{2}$}}
\psfrag{x3}{\small{$x_{3}$}}
\psfrag{ZOOM}{\small{$\text{ZOOM}$}}
\psfrag{Hex}{\small{$\text{Hex}(0)$}}
\psfrag{Hex2}{\small{$\epsilon\text{Hex}(0)$}}
\psfrag{k}{\small{$k_{\epsilon}$}}
\psfrag{k'}{\small{$k_{\epsilon}$}}
\includegraphics[width=9.5cm,height=6.8cm]{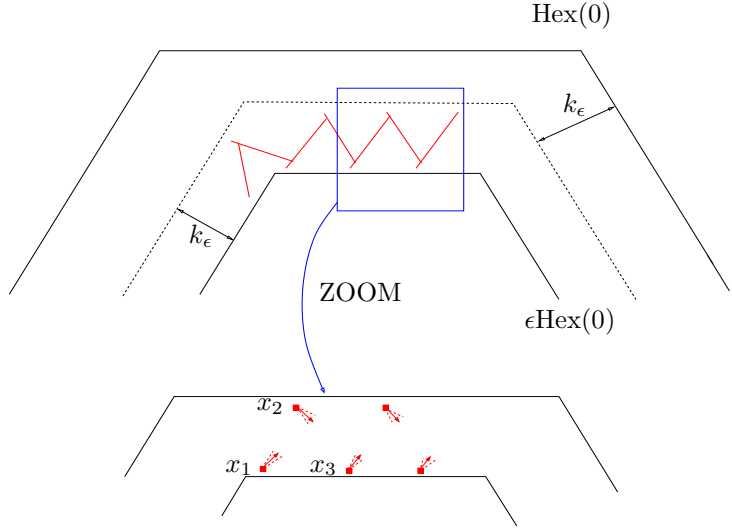}
\caption{\label{fig:pige} {This picture illustrates the building of a loop $\{x_{1},\dots,x_{n}\}$  inside a strip of size $k_{\epsilon}$. Each segment $x_{i}$ of the loop can grow during a time of order $\frac{k_{\epsilon}}{V_{c}(m)}$ without perturbation caused by the segments starting from the outside of $\text{Hex}(0)$}.}
\end{center}
\end{figure}

\begin{dem}{ (Proof of Lemma \ref{lem:lim})}
The proof of Proposition 5.3 in \cite{coupier2016absence} gives
\begin{equation}
\label{eq:moscato}
\mathbf{P} \left( \mathbf{X} \notin E^{(\beta)}_{m} \, | \, 0 \notin \Sigma_{16,m} \right) \le C m^{3} \left( 1 - p_{m}^{6}\right)^{\frac{m}{10}} ~,
\end{equation}
where $C>0$ is a constant. Let us specify that the power $6$ on the probability $p_m$ (defined in (\ref{def:pm})) comes from the fact that we ask the hexagons belonging to two consecutive floors $\text{Cross}_{i}(l)$ and $\text{Cross}_{i+1}(l)$ for some ray $l$ and some index $i\geq i(l)$ (at most $3$ hexagons per floor) to be $(\epsilon,m)$-shield for $\mathbf{X}$. Besides, a geometric construction  leads to a lower bound for $p_m$: for any $m$,
\begin{equation}
\label{eq:minoration}
p_{m} \ge \left( \frac{C_{1}}{V_{c}(m)^{2}} \right)^{\frac{V_{c}(m)}{C_{2}}} ~,
\end{equation}
where $C_{1},C_{2}>0$ are constants. {Let us describe precisely how to obtain \eqref{eq:minoration}. Actually we provide a lower bound for the probability that the  $\text{Graph}_{m}(0)$ creates a barrier with small segments in the ring $\text{Hex}(o)\setminus \epsilon\text{Hex}(o)$ as in the figure \ref{fig:pige}.  We build this barrier  inside the strip $\frac{1+\epsilon}{2}\text{Hex}(0)\setminus\epsilon\text{Hex}(0)$  with thickness $k_\epsilon=\sqrt{3}(1-\epsilon)/2 >0$ and at distance $k_\epsilon$ of  $\text{Hex}(0)^c$. Recall that $\text{Graph}_{m}(0)$ is only composed by grains $\text{Grain}(x,f_{m}(x))$ such that for any $0\le t\le f_{m}(x)\le 1$, the ball $B(H(x,t) , tV_{c}(m)) $ is included $ \text{Hex}(0)$. So the time to grow has to be smaller than $k_\epsilon/V_c(m)$. We fix $\mathscr{V}>0$ such that each segment, with positive probability, has a velocity between  $\mathscr{V}$ and $2\mathscr{V}$. So the length of each segment is of order  $\mathscr{V} k_\epsilon/V_c(m)$ and therefore we need a number of  segments of order $V_c(m)$ for building the barrier. Now we just give the order without the multiplicative constants; the asymptotic is when $m\to\infty$. It remains to localize properly the starting point and the good orientation  of each segment. The location of each starting point can be chosen as a disk with radii of order $1/V_c(m)$ since the length between two consecutive locations is also of order $1/V_c(m)$ (the order of the length of a segment). Each orientation can be chosen randomly inside a deterministic range of angles. We have just to arrange that each segment hits correctly the following segment. Therefore the probability that each segment has a good starting location and orientation is of order  $1/V_c(m)^2$. The independence property of the Poisson point process and the description above guarantee that the probability that  $\epsilon\text{Hex}(0)$ is encircled in a loop is larger than $(C_{1}/V_{c}(m)^{2})^{V_{c}(m)/C_{2}}$, where $C_{1}, C_2$ are positive constant independent of $m$.}

Combining \eqref{eq:moscato} and \eqref{eq:minoration}, we obtain
\begin{equation}
\label{eq:dorian}
\mathbf{P} \left( \mathbf{X}\notin E^{(\beta)}_{m} \, | \, 0 \notin \Sigma_{16,m} \right) \le C m^{3} \exp \Big( -\frac{m}{10} \Big( \frac{C_{1}}{V_{c}(m)^{2}} \Big)^{\frac{6V_{c}(m)}{C_{2}}} \Big) ~.
\end{equation}
Thus, for some $C_3>0$,
$$
\Big( \frac{C_{1}}{V_{c}(m)^{2}} \Big)^{\frac{6V_{c}(m)}{C_{2}}} \geq \exp \big( -C_3 V_{c}(m) \log V_{c}(m) \big)
$$
which is larger than $m^{-1/2}$, for $m$ large enough, using $V_{c}(m)=(\log m^{3})^{1/s}$ with $s>1$. The expected result then follows.
\end{dem}

\subsection{Conclusion}
\label{sect:Step5}

{According to the strategy developped at the end of Section \ref{sect:Step2}, recall that we have assumed that the Forward set $\text{For}(\mathbf{X}_{\gamma},\gamma)$ is infinite with positive probability and it is enough (to get a contradiction) to prove that, with positive probability, $\text{For}(\mathbf{X}_{\gamma},\gamma)$ contains infinitely many $(r,R,W,K,A)$-almost looping points for a deterministic set of parameters $(r,R,W,K,A)$.}

To do it, let us strengthen the event $\mathscr{E}_{m}$ into a new event $\mathscr{E}'_{m}$ with assuming two extra conditions; $(i)$ any almost looping point $x\in\mathbf{X}_{\text{Hex}^{2m}(0)}$ satisfies $\text{radius}(A_{x})>\delta_{m}$ (where $A_x$ is the suitable region to break the Forward set of $x$) and $(ii)$ $\# \mathbf{X}_{\text{Hex}^{4m}(0)} \leq K_{m}$. Like $\mathscr{E}_{m}$, this new event $\mathscr{E}'_{m}$ satisfies the same key properties, i.e. Items $(a)$ \textbf{Localization} and $(b)$ \textbf{ALP property} of Proposition \ref{prop:ALP}, as well as Lemma \ref{lem:lim} provided $\delta_{m}\to 0$ and $K_{m}\to\infty$ fast enough with $m$. Hence, Proposition \ref{prop:NoPerco} still holds with $\mathscr{E}'_{m}$ instead of $\mathscr{E}_{m}$, i.e. $\mathbf{Z}^{2}\!\setminus\!\mathcal{S}_m$ does not percolate for any $m$ large enough, where $\mathcal{S}_m=\{z\in\mathbf{Z}^{2} : z\notin\Sigma_{16,m} \mbox{ and } \tau_{-mz}(\mathbf{X})\in\mathscr{E}_{m}'\}$.

Let us pick such a large $m$. Even if the trajectory $\text{For}(\mathbf{X}_{\gamma},\gamma)$ visits infinitely many connected components of $\mathbf{Z}^{2}\!\setminus\!\mathcal{S}_m$, to go from one of them to another one, $\text{For}(\mathbf{X}_{\gamma},\gamma)$ has to cross a set of blocks on which $\mathscr{E}'_{m}$ occurs and then admits inside {a $(5m,6m,V_{c}(m),K_y,A_y)$-almost looping point $y$ by the \textbf{ALP property}. Thus, combining with condition $(ii)$ of the strengthened event $\mathscr{E}_{m}'$, we obtain that with positive probability, $\text{For}(\mathbf{X}_{\gamma},\gamma)$ contains infinitely many $(5m,6m,V_{c}(m),K_m,A_y)$-almost looping points $y$. Note that only the fifth parameter $A_y$ may depend on $y$ and be random.}

{To make deterministic this fifth parameter, we use condition $(i)$ of $\mathscr{E}_{m}'$.} Let $y\in\text{For}(\mathbf{X}_{\gamma},\gamma)$ be such a visited almost looping point. There exists $z_{y}\in\mathbf{Z}^{2}$ such that $\tau_{-mz_{y}}(\mathbf{X})\in\mathscr{E}'_{m}$ and $[y,h_{g}(\mathbf{X}_{\gamma},y)]\subset\text{Hex}^{2m}(mz_{y})$. As previously mentioned, $y$ is a $(5m,6m,V_{c}(m),K_{m},A_y)$-almost looping point where the random ball $A_{y}$, with radius larger than $\delta_m$, is included in $(B(\eta,5m)\times[0,2\pi]\times[0,V_{c}(m)])^{3}$. Let us consider a finite covering of $(B(0,5m)\times[0,2\pi]\times[0,V_{c}(m)])^{3}$ by open Euclidean balls $\{\mathscr{K}_{j}, 1\le j\le j(m)\}$ with radius $\frac{\delta_{m}}{2}$. Then, the pigeonhole principle asserts that there exists a deterministic ball $\mathscr{K}_{j_{0}}$ with some deterministic $1\le j_0\le j(m)$ such that, with positive probability, among the $(5m,6m,V_{c}(m),K_{m},A_y)$-almost looping points visited by $y\in\text{For}(\mathbf{X}_{\gamma},\gamma)$, infinitely many of them satisfy $\tau_{y}(\mathscr{K}_{j_{0}}) \subset A_{y}$. We then conclude that with positive probability, $\text{For}(\mathbf{X}_{\gamma},\gamma)$ contains infinitely many $(5m,6m,V_{c}(m),K_{m},\mathscr{K}_{j_{0}})$-almost looping points of $\mathbf{X}_{\gamma}$.

\section*{Acknowledgement}

The authors thank J.-B. Gou\'er\'e for fruitful discussions on this problem. This work was supported in part by the Labex CEMPI (ANR-11-LABX-0007-01), the CNRS GdR 3477 GeoSto and by the ANR project PPPP(ANR-16-CE40-0016).

\addcontentsline{toc}{section}{Bibliography}

\bibliographystyle{plain}
\bibliography{biblio2}

\begin{thebibliography}{1}

\bibitem{athreya2004branching}
K.~B. Athreya and P.~E. Ney.
\newblock {\em Branching processes}.
\newblock Courier Corporation, 2004.

\bibitem{chiu2013stochastic}
S.~N. Chiu, D.~Stoyan, W.~S. Kendall, and J.~Mecke.
\newblock {\em Stochastic geometry and its applications}.
\newblock John Wiley \& Sons, 2013.

\bibitem{coupier2016absence}
D.~Coupier, D.~Dereudre, and S.~Le~Stum.
\newblock Absence of percolation for poisson outdegree-one graphs.
\newblock {\em Annales de l'Institut Henri Poincar{\'e}, Probabilit{\'e}s et
  Statistiques}, 56(2):1179--1202, 2020.

\bibitem{daley2014two}
D.~Daley, S.~Ebert, and G.~Last.
\newblock Two lilypond systems of finite line-segments.
\newblock {\em Probability and Mathematical Statistics}, 36(2):221--246, 2016.

\bibitem{daley2005descending}
D.~Daley and G.~Last.
\newblock Descending chains, the lilypond model, and mutual-nearest-neighbour
  matching.
\newblock {\em Advances in applied probability}, pages 604--628, 2005.

\bibitem{gouere2008subcritical}
J.-B. Gou{\'e}r{\'e}.
\newblock Subcritical regimes in the poisson boolean model of continuum
  percolation.
\newblock {\em The Annals of Probability}, pages 1209--1220, 2008.

\bibitem{hall1985continuum}
P.~Hall.
\newblock On continuum percolation.
\newblock {\em The Annals of Probability}, pages 1250--1266, 1985.

\bibitem{liggett1997domination}
T.~Liggett, R.~Schonmann, and A.~Stacey.
\newblock Domination by product measures.
\newblock {\em The Annals of Probability}, 25(1):71--95, 1997.

\bibitem{meester1996continuum}
R.~Meester and R.~Roy.
\newblock {\em Continuum percolation}, volume 119.
\newblock Cambridge University Press, 1996.

\end{thebibliography}
\end{document}